%% file: mollified-coll-wiley.tex
\documentclass[AMA,Times1COL]{WileyNJDv5} %STIX1COL,STIX2COL,STIXSMALL

\usepackage{amsmath}
\usepackage{amssymb}
\usepackage{mathtools}
\usepackage[utf8]{inputenc}
\usepackage{subfig}

\usepackage{csml}

\articletype{}
\startpage{1}
\begin{document}

\title{Point collocation with mollified piecewise polynomial approximants for high-order partial differential equations}

\author[{1,2}]{Dewangga Alfarisy}
\author[1]{Lavi Zuhal}
\author[3]{Michael Ortiz}
\author[4]{Fehmi Cirak}
\author[{1,2}]{Eky Febrianto}

\authormark{ALFARISY \textsc{et al.}}
\titlemark{Point collocation with mollified piecewise polynomial approximants for high-order partial differential equations}

\address[1]{Faculty of Mechanical and Aerospace Engineering, Bandung Institute of Technology, Bandung, 40116, Indonesia}
\address[2]{Glasgow Computational Engineering Centre, University of Glasgow, Glasgow, G12 8QQ, UK}
\address[3]{Graduate Aerospace Laboratories, California Institute of Technology, Pasadena, CA,  91125, USA}
\address[4]{Department of Engineering, University of Cambridge, Cambridge, CB2 1PZ, UK}

\corres{Corresponding author: Eky Febrianto, \email{eky.febrianto@glasgow.ac.uk}}

\abstract[Abstract]{The solution approximation for partial differential equations (PDEs) can be substantially improved using smooth basis functions. The recently introduced mollified basis functions are constructed through mollification, or convolution, of cell-wise defined piecewise polynomials with a smooth mollifier of certain characteristics. The properties of the mollified basis functions are governed by the order of the piecewise functions and the smoothness of the mollifier. In this work, we exploit the high-order and high-smoothness properties of the molli- fied basis functions for solving PDEs through the point collocation method. The basis functions are evaluated at a set of collocation points in the domain. In addi- tion, boundary conditions are imposed at a set of boundary collocation points distributed over the domain boundaries. To ensure the stability of the resulting linear system of equations, the number of collocation points is set larger than the total number of basis functions. The resulting linear system is overdeter- mined and is solved using the least square technique. The presented numerical examples confirm the convergence of the proposed approximation scheme for Poisson, linear elasticity, and biharmonic problems. We study in particular the influence of the mollifier and the spatial distribution of the collocation points. 
}

\keywords{collocation, convolution, high-order, mollification, polytopic meshes, smooth basis functions}

\maketitle

%\renewcommand\thefootnote{}
%\footnotetext{\textbf{Abbreviations:} PDE, partial differential equation; PCM, point collocation method.}

%\renewcommand\thefootnote{\fnsymbol{footnote}}
%\setcounter{footnote}{1}

\input{introduction}

\input{mollified}

\input{collocation}

\input{examples}

\input{conclusions}

%\backmatter
%\bmsection*{Author contributions}

%\bmsection*{Acknowledgments}
%\bmsection*{Financial disclosure}
%None reported.

\bmsection*{Conflict of interest statement}

The authors declare no potential conflict of interests.

\bmsection*{Data availability statement}

The data that support the findings of this study are available from the corresponding author upon reasonable request.

\bibliography{mollified-coll}

%\bmsection*{Supporting information}
%Additional supporting information may be found in the online version of the article at the publisher’s website.

% \appendix

\end{document}

%% file: introduction.tex
%--------------------------------------------------------------------------------
\section{Introduction}
\label{sec:introduction}
%--------------------------------------------------------------------------------
%

% general of solving PDE using smooth basis functions - through FEM and IGA
\subsection{Motivation}
\label{sec:motivation}

Finding approximate solutions to high-order partial differential equations (PDEs) is a foundational task in various scientific and engineering problems, including gradient theories of elasticity and plasticity~\cite{rudraraju2014three, codony2019immersed, khakalo2022strain}, phase-field modelling of sharp interfaces~\cite{gomez2008isogeometric, liu2013isogeometric, ambati2015review}, and plate and shell models~\cite{Cirak2000, kiendl2009isogeometric, cirak2014computational}. The strong form of these PDEs typically impose stringent smoothness requirements on the approximation schemes. In certain scenarios, it is advantageous to pursue solutions that adhere to the weak form of PDEs, where the smoothness requirements are less stringent. Nevertheless, constructing an approximate solution for these PDEs commonly involves an ascending sequence of smooth basis functions, which are a subset of the solution space.

On a parallel note, there has been a growing interest in isogeometric analysis (IGA) that employs smooth basis functions prevalent in computer-aided design (CAD), such as the spline-based NURBS~\cite{piegl1997nurbs, farin2001} and subdivision surfaces~\cite{Peters2008}, to solve various PDEs. These smooth basis functions found applications within the framework of the finite element method (FEM)~\cite{Cirak2000, Cirak:2002aa, Hollig2003, Hughes2005, Kamensky2015a} based on the weak form of PDEs and the Galerkin method. Although the primary aim of IGA is to streamline design and analysis, it has undoubtedly leveraged the broader utilisation of smooth basis functions in computational mechanics. 

Solving the strong form of PDEs through the point collocation method (PCM) is feasible when using sufficiently smooth basis functions. In PCM, the PDEs are directly evaluated at specific spatial points, known as collocation points, with boundary conditions imposed at the boundary collocation points. Because no integration is required for evaluating the strong form, collocation presents a straightforward and usually cost-effective alternative to the traditional FEM~\cite{Schillinger2013, Hillman2018}.

% collocation with smooth basis functions - a short review
\subsection{Previous work}
\label{sec:prev-work}

Various smooth basis functions have been employed in the collocation framework, and many of them are also known in the context of IGA, such as B-splines~\cite{Manni2015, Torre2023immersed, Torre2023isogeometric}, T-splines~\cite{Casquero2016}, and NURBS~\cite{Auricchio2010, Auricchio2012,  Anitescu2015}. These functions are typically mesh-based, and their smoothness may degrade when the tensor product structure of the mesh breaks, particularly at extraordinary vertices and edges~\cite{farin2001, Toshniwal2017, zhang2018subdivision, ZHANG2020112659, koh2022optimally, verhelst2024comparison}. Mesh-free approximants, including moving least squares~\cite{Wu2013}, radial basis functions~\cite{Zhang2000, hu2007}, reproducing kernel~\cite{Aluru2000, Chi2013}, and the maximum entropy (max-ent) approximants~\cite{Fan2018, greco2020, Fan2021} have also been applied for point collocation. Most of these approximants are not polynomial and often need to satisfy certain consistency criteria~\cite{liu1995reproducing, Bonet2000} to maintain polynomial reproducibility. 

% mollified approximants
Smooth basis functions can be easily constructed through mollification, or convolution, of piecewise polynomials with a smoothing kernel~\cite{adams2003sobolev}. In a conceptual sense, mollification is fundamentally different from traditional interior approximation schemes. Interior methods, such as FEM, approximate the solution using interpolants belonging to the solution space, for example~$H^1$. In contrast, exterior methods, such as mollification, enable solution approximation using functions outside of the solution space, that is, using a piecewise polynomial basis. While interior methods are well-established, exterior methods, such as the mollification approach, are simpler, more general, and relatively newer. As highlighted in Febrianto et.al.~\cite{Febrianto2021}, mollified approximants maintain the order of the piecewise polynomial basis, while improving smoothness as determined by the kernel, or mollifier. This gives rise to an appealing characteristic of the mollified basis functions where the smoothness and polynomial order can be arbitrary. This smooth basis construction strategy shares similarities with the convolutional definition of B-splines~\cite{de1986b, Sabin:2010aa} and those of simplex splines~\cite{micchelli1995, grandine1988stable}. 

In the mollified approach, the polynomials are defined over cells or meshes, enabling a faithful evaluation of the convolution integral, particularly when employing a compactly supported polynomial mollifier. The resulting basis functions span the same space as the piecewise polynomial, therefore avoiding the need for corrections that imposes the polynomial reproducibility on the kernel~\cite{liu1995reproducing, Bonet2000, Liu2004, Li2004}. Furthermore, unlike in most mesh-based approximants, the construction of mollified basis functions is not restricted to a specific type of domain discretisation, and their smoothness remains unaffected by extraordinary vertices. This versatility is particularly advantageous when a specific type of discretisation is faster and more robust to obtain, such as the Cartesian grid or Voronoi tessellation~\cite{Du1999, Ray2018, Abdelkader2020, Bishop2020, Sukumar2022}, which arguably aligns with the objective of reducing the cost of domain discretisation. 

% work description - collocation method using mollified approximant
\subsection{Contributions}
\label{sec:contribution}

\begin{figure}[]
	\centering
	\subfloat[][Domain definition \label{fig:dodMeshCoord}]{
		\includegraphics[scale=0.05]{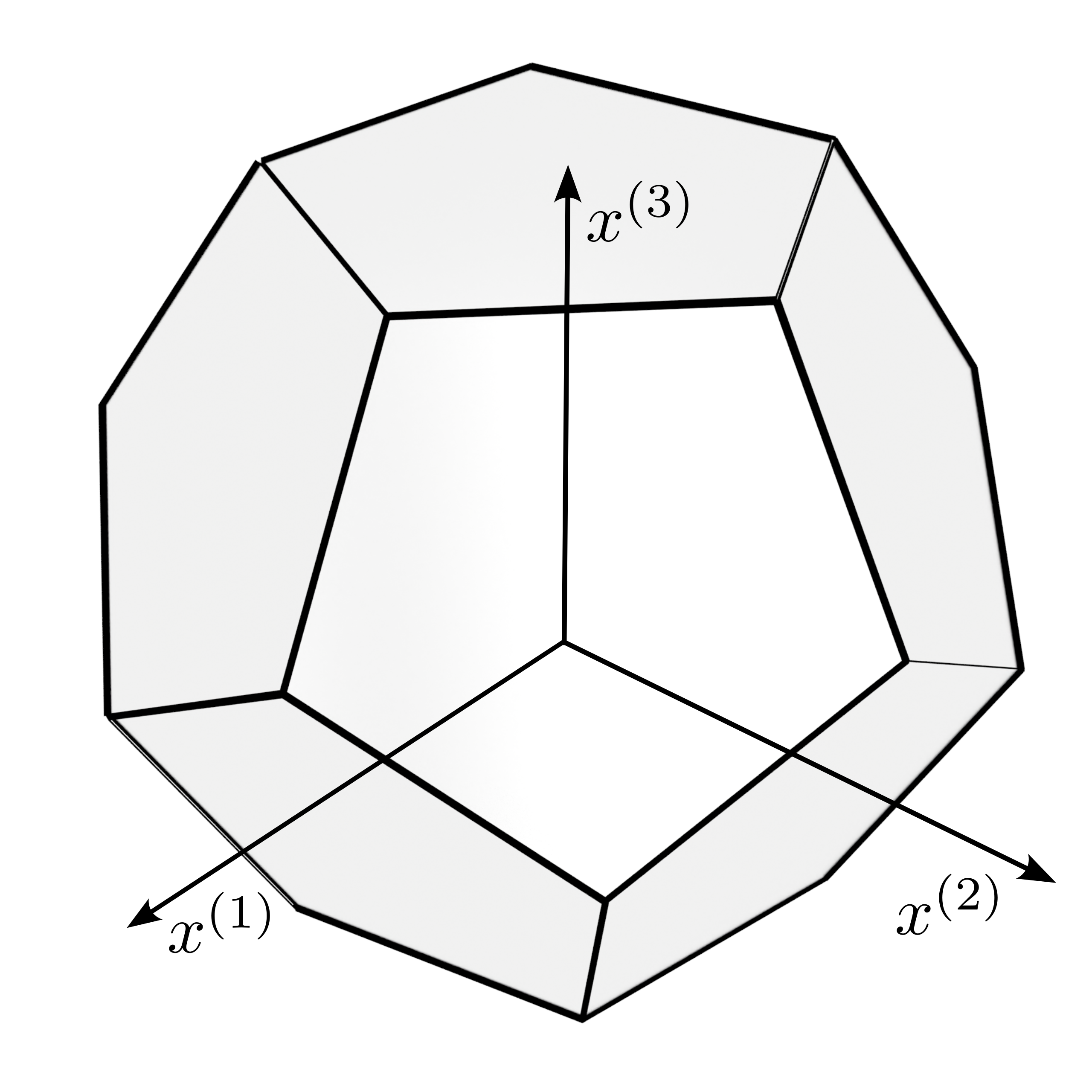}
	}
	\hspace{0.05\textwidth}
	\subfloat[][Voronoi tessellation \label{fig:dodMesh}]{
		\includegraphics[scale=0.07]{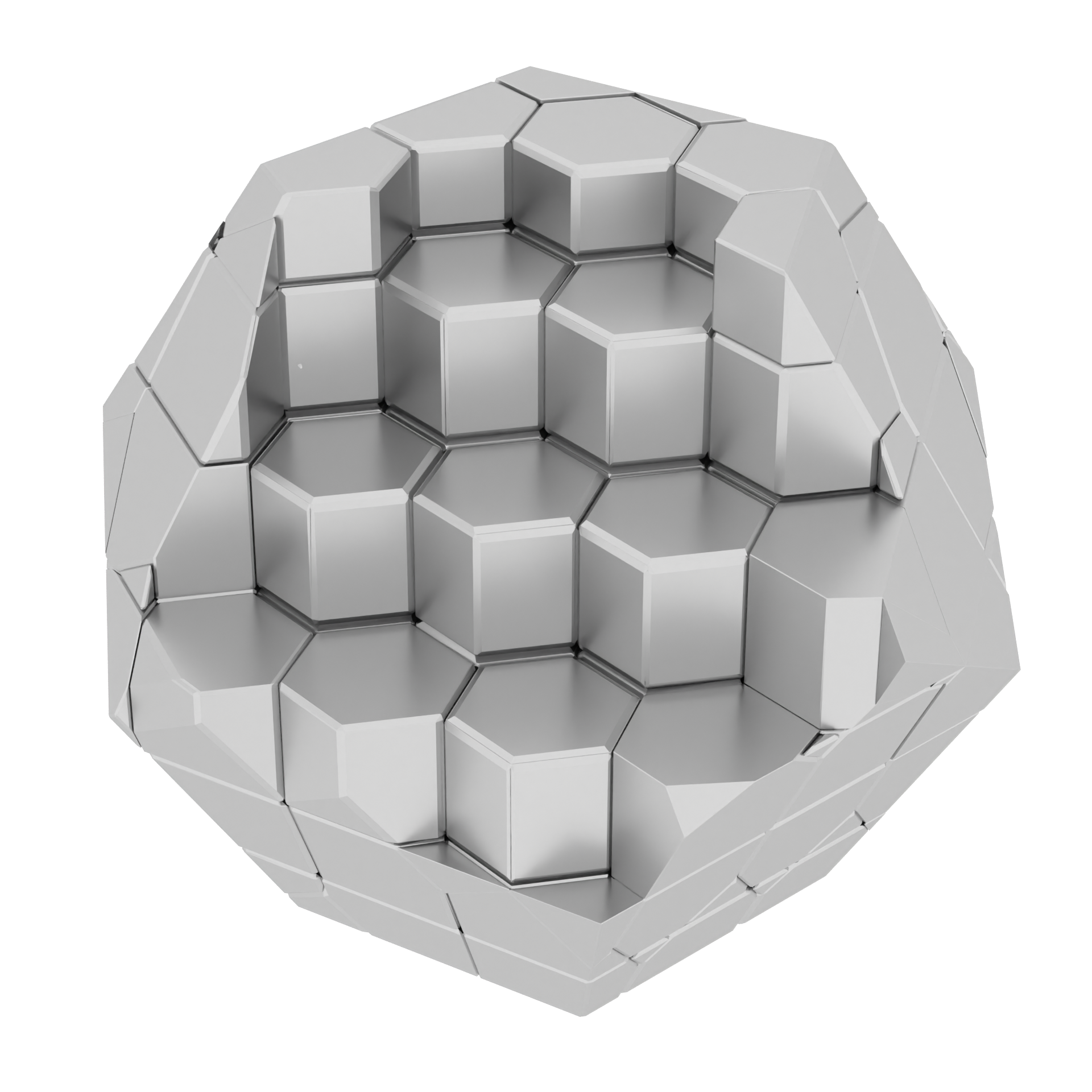}
	}
	\hspace{0.05\textwidth}
	\subfloat[][Solution contour \label{fig:dodResult}]{
		\includegraphics[scale=0.07]{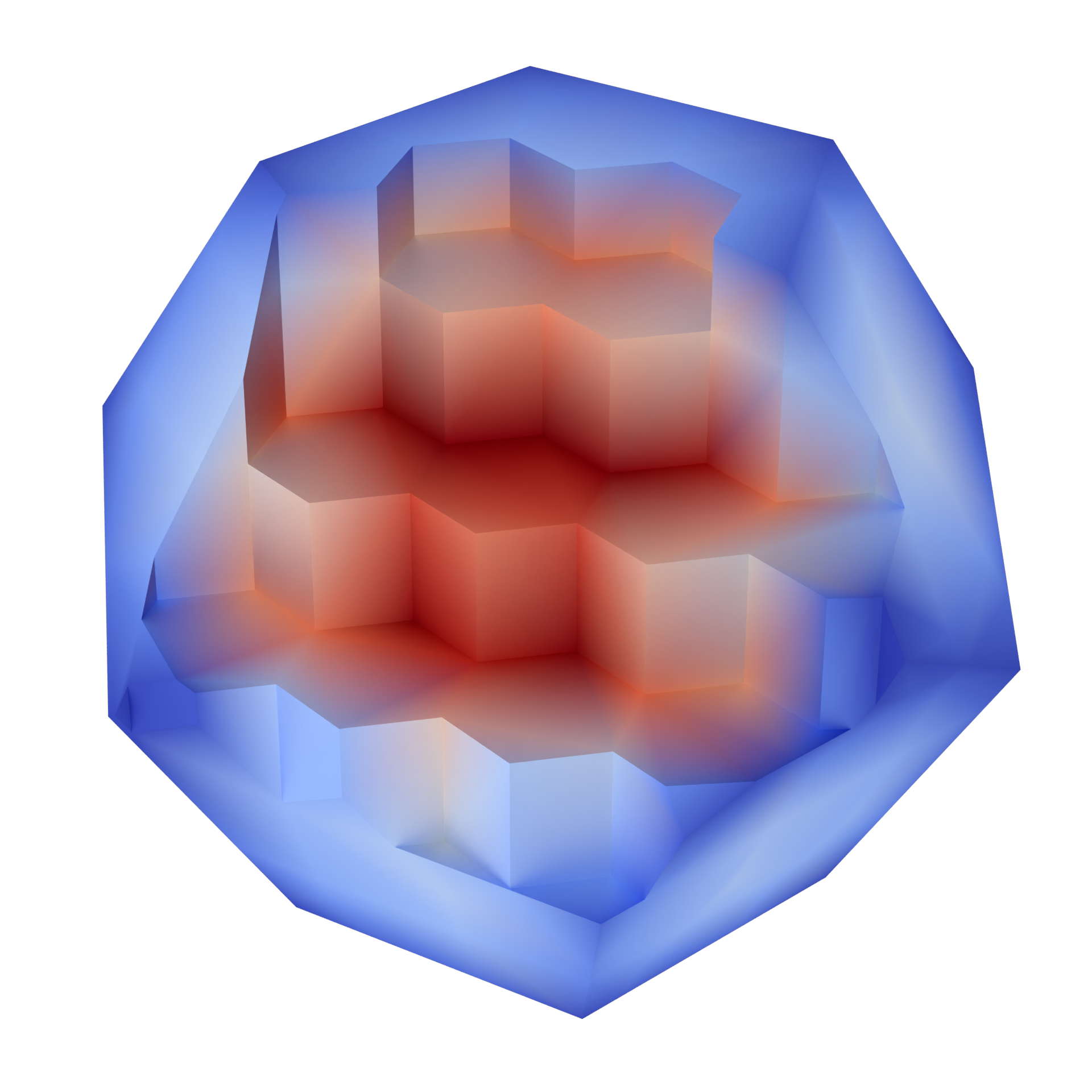}
	}
	\caption[]{The proposed mollified-collocation method. First the domain (a) is discretised into a set of polytopic cells through Voronoi tessellation (b). Each cell has an associated piecewise linear polynomial and a $C^{k-1}$ smooth kernel is used to obtain smooth high-order basis functions. The resulting basis functions are used to solve a $k$-th order PDE problem giving the solution (c).}
	\label{fig:dod}
\end{figure}

In this paper, we present a point collocation method for solving the strong form of PDEs using mollified basis functions. The high smoothness of the mollified basis functions makes them particularly suitable for collocation. Compared to its Galerkin implementation~\cite{Febrianto2021}, the mollified-collocation proposed in this work is simpler and more straightforward to implement. In particular, the mollified-collocation sidesteps the need for accurate and variationally consistent integrations when evaluating the domain and surface integrals in the weak form~\cite{Chen2013, hillman2015}. Furthermore, the mollified-collocation facilitates the strong imposition of boundary conditions, thus circumventing the need for auxiliary methods required to impose stably Dirichlet boundary conditions~\cite{Nitsche1971, Ruberg:2011aa,  boiveau2015penalty, Febrianto2021}. In this work, we numerically investigate the convergence of the approximate solutions of the Poisson, linear elasticity, and biharmonic problems over polytopic elements, see for example Figure~\ref{fig:dod}. Additionally, we assess the effect of the polynomial order, as well as mollifier smoothness and width on the accuracy of numerical approximation.  

The unstructured shape of the support of mollified basis functions is incompatible with knot-based collocation point distributions, such as the Greville and Demko abscissae, commonly employed in IGA collocation methods~\cite{Johnson2005, Auricchio2010, Manni2015, Wang2021}. Moreover, in the proposed mollified-collocation, multiple basis functions might overlap over each cell, implying that collocating solely at a cell's centroid or the Voronoi seeds would result in an underdetermined system matrix. To address this challenge, more than one collocation point in a cell is distributed according to the selected scheme. Additionally, a constraint is imposed to ensure that the total number of collocation points exceeds the number of basis functions involved in the computation, a quantity that can be predetermined based on the number of cells and the polynomial order. This approach results in an overdetermined (non-square) system matrix, which can be solved using a standard least square technique. In this study, we explore three schemes for distributing the collocation points: uniform, Gauss quadrature, and quasi-random schemes. Furthermore, we analyse their effects on the convergence of the approximation error. Specifically, for quasi-random collocation points, we perform stochastic experiments, and sample perturbations $n$ times to ascertain the mean and standard deviation of errors.

\subsection{Overview}
\label{sec:overview}

The structure of this paper is as follows. First, we revisit the mollification approach for smoothing piecewise polynomials. We then apply this principle to construct smooth basis functions in one and higher dimensions. Subsequently, we discuss the evaluation of the basis functions and their derivatives at a point in space. Later, we describe the use of the mollified basis functions in the collocation framework, including considerations on distributing collocation points in space. Finally, we present several numerical Poisson, linear elasticity, and biharmonic examples in one, two, and three dimensions. 

%% file: mollified.tex
%--------------------------------------------------------------------------------
\section{Review of mollified piecewise polynomial approximants}
\label{sec:mollified}
%--------------------------------------------------------------------------------
%
In this section, we review the mollified piecewise polynomial approximants used to discretise PDEs. We begin by describing the mollification of piecewise polynomial functions, resulting in a global function smoother than the smoothing kernel, that is, mollifier. This characteristic is particularly advantageous when approximating the solution of PDEs using the collocation method. Next, we derive the mollified basis functions based on the piecewise polynomial defined in each domain partition, referred to as a cell. This work focuses on polytopic meshes, such as the Voronoi tessellation, for domain partitioning. The values of the basis functions can be obtained at a point in space by evaluating a convolution integral, which is also detailed in this section. 
%
%--------------------------------------------------------------------------------
\subsection{Mollification of the piecewise polynomial}
\label{sec:mollification}
%--------------------------------------------------------------------------------
%
For brevity, we illustrate the mollification of piecewise polynomials in a one-dimensional setting. We begin by considering the domain~$\Omega \in \mathbb{R}^1$ discretised into a set of~$n_c$ non-overlapping cells~$\{ \omega_i \}$, such that,
\begin{equation}\label{eq:partitioning}
	\Omega = \bigcup_{i=1}^{n_{c}} \omega_i  \, .
\end{equation}
On each cell~$\omega_i$, we define a local polynomial 
\begin{equation} \label{eq:fi}
	f_i(x) = 
	\begin{cases}
		\vec p_i (x)  \cdot \vec \alpha_i \quad  & \text{if } \,  x \in \omega_i \\
		0						 &  \text{if } \,  x \not \in \omega_i 
	\end{cases} \, , 
\end{equation}
where~$\vec p_i(x)$ is a vector containing a local polynomial basis of order~$r_p$ and~$\vec \alpha_i$ are the respective coefficients of the basis. Common choices for the basis $\vec p_i(x)$ include, but are not limited to, monomials, Lagrange polynomials, and Bernstein polynomials. While it is possible to vary the polynomial order in each cell, this study exclusively focuses on a uniform polynomial order for all cells. The sum of the local polynomials defined over the entire domain~$\Omega$ constitutes the global polynomial 
\begin{equation} \label{eq:fx}
	f(x) = \sum_i \vec p_i(x) \cdot \vec \alpha_i \, . 
\end{equation}
Note that across the cell boundaries this function will be discontinuous. 

We consider the smoothing of the piecewise polynomial~$f(x)$ through convolution, referred to as \textit{mollification}, with a smooth kernel referred to as a \textit{mollifier}. The mollification of  $f(x)$ with a mollifier~$m(x)$ is defined as
\begin{equation} \label{eq:mollify}
	\widehat f (x) = m(x) * f(x)  = \int_\Omega m(x- y) \, f(y) \D y   \, .
\end{equation}
We require the mollifier to be non-negative, have a unit volume, and have finite support per our previous work~\cite{Febrianto2021}. An important characteristic of mollification is that the mollified functions $\widehat f (x)$ can exactly reproduce polynomials of order~$r_p$~\cite{Febrianto2021}. 

When the derivative of the mollifier~$m(x)$ exists, the $k$-order derivative of the mollified function~$\widehat{f}(x)$ is given by
\begin{equation} \label{eq:mollifyDeriv}
	\frac{\D^{k}}{\D^{k} x} \widehat{f}(x) =  \int_\Omega \frac{\D^{k} m(x - y)}{\D^{k} x} \,  f(y) \D y \, .
\end{equation}
For the mollified function $\widehat{f}(x)$ to be $C^{k}$-smooth, suitable for instances where it is considered a candidate solution for a $k$-th order PDE, the mollifier must be $C^{k-1}$-smooth. Note that the original function $f(x)$ is discontinuous across the cell boundaries. 

\begin{figure}[tb]
	\centering
	\subfloat[][Mollification with a symmetric kernel] {
		\includegraphics[width=0.46 \textwidth]{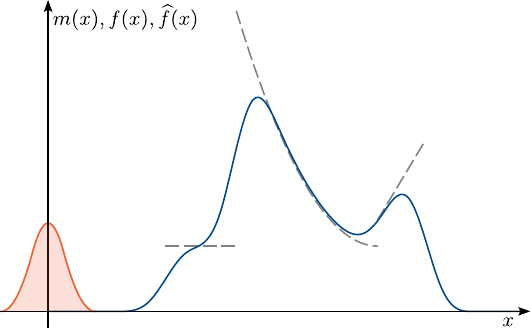} \label{fig:gradPlot}} 
	\hfill
	\subfloat[][Mollification with an asymmetric kernel] {
		\includegraphics[width=0.46 \textwidth]{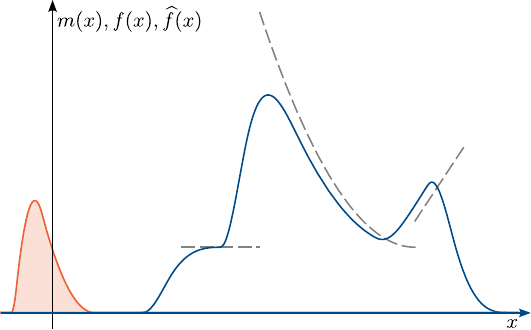} \label{fig:hessPlot}} 
	\caption[Mollification of piecewise linear functions]{Mollification of piecewise linear functions with $C^1$-smooth quadratic B-spline mollifier. The resulting function $\widehat{f}(x)$ is $C^2$-smooth. Mollification with symmetric and asymmetric kernels is used to exemplify the parallel between the convolution operator and feature maps in CNN. \label{fig:convolution}}
\end{figure}

An intriguing observation is the parallel between mollification~\eqref{eq:mollify} and the widely used convolutional neural networks (CNN). A single CNN layer can be represented as $\widehat f (x) = \sigma (\sum_j m(x - x_j) f(x_j)  + b )$, where $\sigma$ represents the nonlinear activation function and $b$ is the bias. This equation can be viewed as the discrete form of the mollification in~\eqref{eq:mollify}, specifically,~$\widehat f (x) = \sum_j m(x - x_j) f(x_j) w_j$, where $x_j$ and $w_j$ are the position and weight of the quadrature points, respectively. In CNNs, the convolution operator is often referred to as a \textit{feature map}, with the convolution kernel containing shared weights, or hyperparameters, which require learning. As an example,  Figure~\ref{fig:convolution} illustrates the mollification of a piecewise polynomial function $f(x)$ with a $C^1$-smooth symmetric (left) and asymmetric mollifier (right). Moreover, choosing a different mollifier here is equivalent to selecting a different feature map in CNN. Although the parallels are evident, this paper does not pursue this comparison further. 

%--------------------------------------------------------------------------------
\subsection{Mollified basis functions}
\label{sec:mollified-basis}
%--------------------------------------------------------------------------------
%
We now use the mollification approach to derive uni- and multivariate basis functions. In the univariate case, following the previous discussion in Section~\ref{sec:mollification}, we partition the domain~$\Omega$ into a set of non-overlapping cells~$\{\omega_i \}$ where the piecewise polynomials $\{ f_i(x) \}$ are defined. By introducing the piecewise definition of the global polynomial~$f(x)$ from~\eqref{eq:fx} into the definition of mollification~\eqref{eq:mollify}, we obtain
\begin{equation} \label{eq:mollifyCells}
	\widehat{f}(x) = \sum_i \vec \alpha_i \cdot \int_{\omega_i} m(x - y) \, \vec p_i(y) \D y \, .
\end{equation}
Here, a polynomial~$ f_i(x)$ is zero outside of the respective cell~$\omega_i$. This allows us to consider as an integration domain~$\omega_i$ instead of the whole domain $\Omega$. We can then express~\eqref{eq:mollifyCells} as a linear combination of basis functions and their coefficients
\begin{equation} 
	\widehat{f}(x) = \sum_i \vec \alpha_i \cdot \vec N_i(x) \, .
\end{equation}
The mollified basis function $\vec N_i(x)$ is obtained by convolving the piecewise polynomial basis $\vec p_i (y)$ belonging to each cell $\omega_i$ with the mollifier $m(x)$ such that
\begin{equation} \label{eq:basis}
	\vec N_i(x) =  \int_{\omega_i} m(x - y) \, \vec p_i (y)  \D y  \, ,
\end{equation}
where the convolution is individually evaluated for each component of $\vec p_i (x)$. In this paper, we consider the vector $\vec p_i(x)$ as the monomial basis defined locally in each cell $\omega_i$. The local monomials are centred at the centroid $c_i$ of each cell, that is,  
\begin{equation} \label{eq:scaledMono}
	\vec p_i(x) = \begin{pmatrix}
		1 & \xi & \xi^2  & \xi^3  & \ldots 
	\end{pmatrix}
	\quad \text{  with } \xi =  \frac{2(x- c_i)}{h_c} \, ,
\end{equation}
where~$h_c$ is the cell size. The scaling by $2 / h_c$ ensures that all mollified basis functions have a similar maximum value, which improves the conditioning of the system matrix. The derivatives of the basis functions can be obtained using~\eqref{eq:mollifyDeriv} by considering the derivative of the mollifier
\begin{equation} \label{eq:basisDeriv}
	\frac{\D}{\D x} \vec N_i (x) =  \int_\Omega \frac{\D m(x - y)}{\D x} \,  \vec p_i (y) \D y \, .
\end{equation}

\begin{figure}[] 
	\centering
	\subfloat[$N^0_i(x)$]{\includegraphics[width=0.38\textwidth]{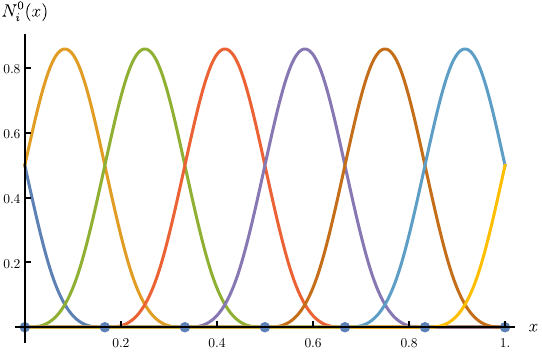} \label{fig:basisAll0B2Unif}}
	\hspace{0.1\textwidth}
	\subfloat[$\frac{\D^2 N^0_i(x)}{\D x^2}$] {\includegraphics[width=0.38\textwidth]{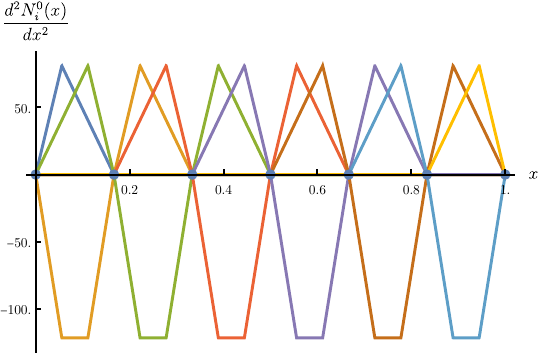} \label{fig:hessAll0B2Unif}} \\
	\subfloat[$N^1_i(x)$]{\includegraphics[width=0.38\textwidth]{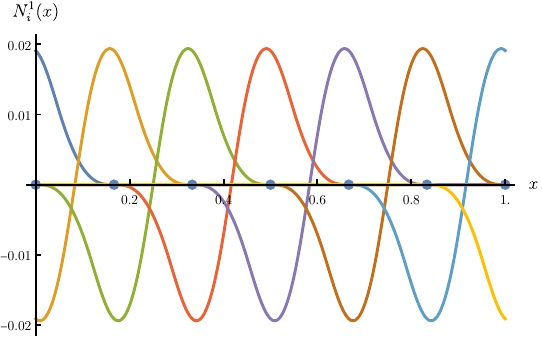} \label{fig:basisAll1B2Unif}} 
	\hspace{0.1\textwidth}     
	\subfloat[$\frac{\D^2 N^1_i(x)}{\D x^2}$] {\includegraphics[width=0.38\textwidth]{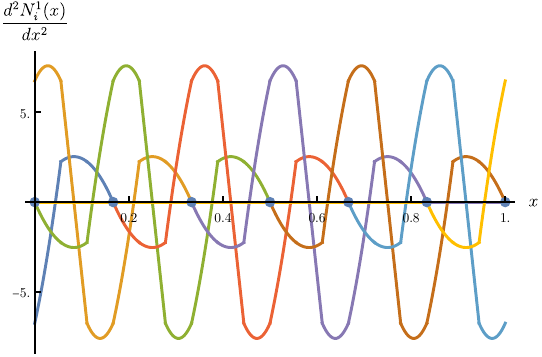} \label{fig:hessAll1B2Unif}} \\
	\subfloat[$N^2_i(x)$]{\includegraphics[width=0.38\textwidth]{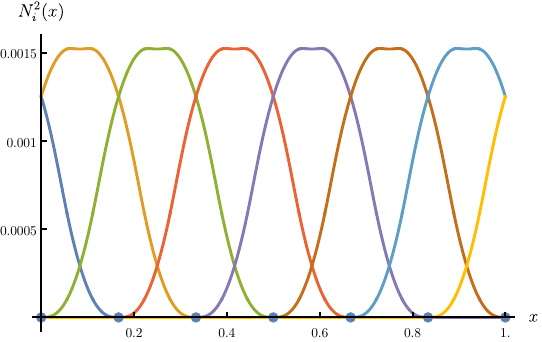} \label{fig:basisAll2B2Unif} }
	\hspace{0.1\textwidth} 
	\subfloat[$\frac{\D^2 N^2_i(x)}{\D x^2}$] {\includegraphics[width=0.38\textwidth]{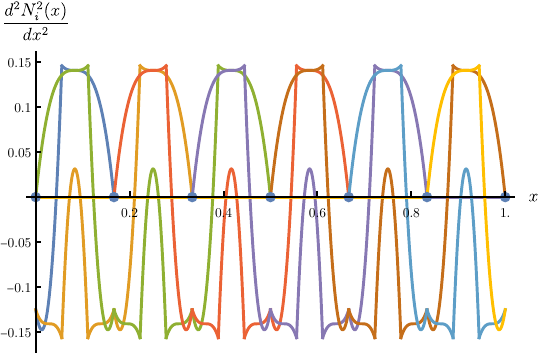} \label{fig:hessAll2B2Unif}} \\
	\subfloat[$N^3_i(x)$]{\includegraphics[width=0.38\textwidth]{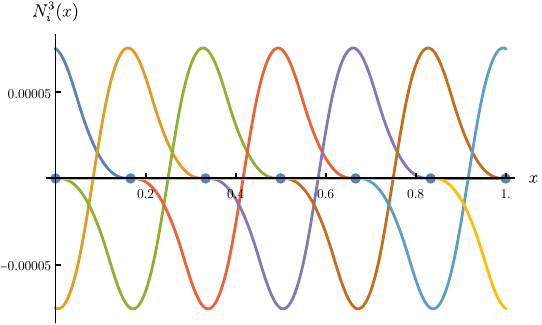} \label{fig:basisAll3B2Unif}}
	\hspace{0.1\textwidth}       
	\subfloat[$\frac{\D^2 N^3_i(x)}{\D x^2}$] {\includegraphics[width=0.38\textwidth]{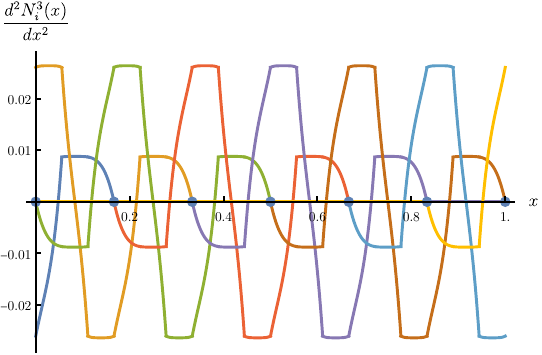} \label{fig:hessAll3B2Unif}}
	\caption[Basis functions and second derivatives using quadratic B-spline kernel]{Basis functions (left column) over $n_c = 6$ uniform cells for $r_p \in \{ 0, 1, 2, 3 \}$ using a quadratic B-spline mollifier and their corresponding second derivatives (right column). The blue dots situated along the $x$-axis denote the cell boundaries. \label{fig:basisAllUnif}}
\end{figure}

As an illustrative example, we consider a one-dimensional domain $\Omega = (0, \, 1)$ discretised into $n_c = 6$ non-overlapping cells with uniform spacing $h_c = 1/6$. Piecewise polynomials up to degree $r_p = 3$ consisting of monomial basis functions defined over each cell are mollified with a quadratic B-spline mollifier with a width $h_m = 2 \, h_c$ and a unit volume. Figure~\ref{fig:basisAllUnif} depicts the obtained mollified basis functions and their respective second derivatives. We stress that each mollified basis can reproduce polynomials up to degree $r_p$. Collocating at positions where either the basis or the second derivative is zero is often undesirable. For this specific example, such zero-valued positions can be predicted, as shown in Figure~\ref{fig:basisAllUnif}. For instance, the linear $N_i^1(x)$ and cubic $N_i^3(x)$ are zero at the cell's centre, and the second derivative of $N_i^0(x)$ coincides with the cell boundary. However, predicting such locations becomes challenging for non-uniform arrangements in both univariate and multivariate cases, especially for polytopic partitions.

Without loss of generality, in the following we continue our discussion on multivariate basis functions, focusing on the bivariate case. In the two-dimensional setting, for the mollifier, we consider the tensor-product of its one-dimensional description
\begin{equation}
	m(\vec x) = m \left(x^{(1)} \right)  \cdot m \left(x^{(2)} \right) \, .
\end{equation}
Consequently, the support of the mollifier~$\Box_{\vec  x}$ is a square, and generally, a hypercube in the multi-dimensional case. Similar to the univariate case, the domain $\Omega$ is subdivided into non-overlapping cells $\{\omega_i\}$. For the multivariate case considered in this work, a polytopic discretisation of the domain $\Omega$ is used. Nevertheless, other non-overlapping partitions, such as Delaunay triangulation/tetrahedralisation, quadrilateral/hexahedral mesh, and Cartesian grid, are also viable options. 

\begin{figure}[tb]
	\centering
	\subfloat[][Mollifier and constant local basis functions] {
		\includegraphics[width=0.35 \textwidth]{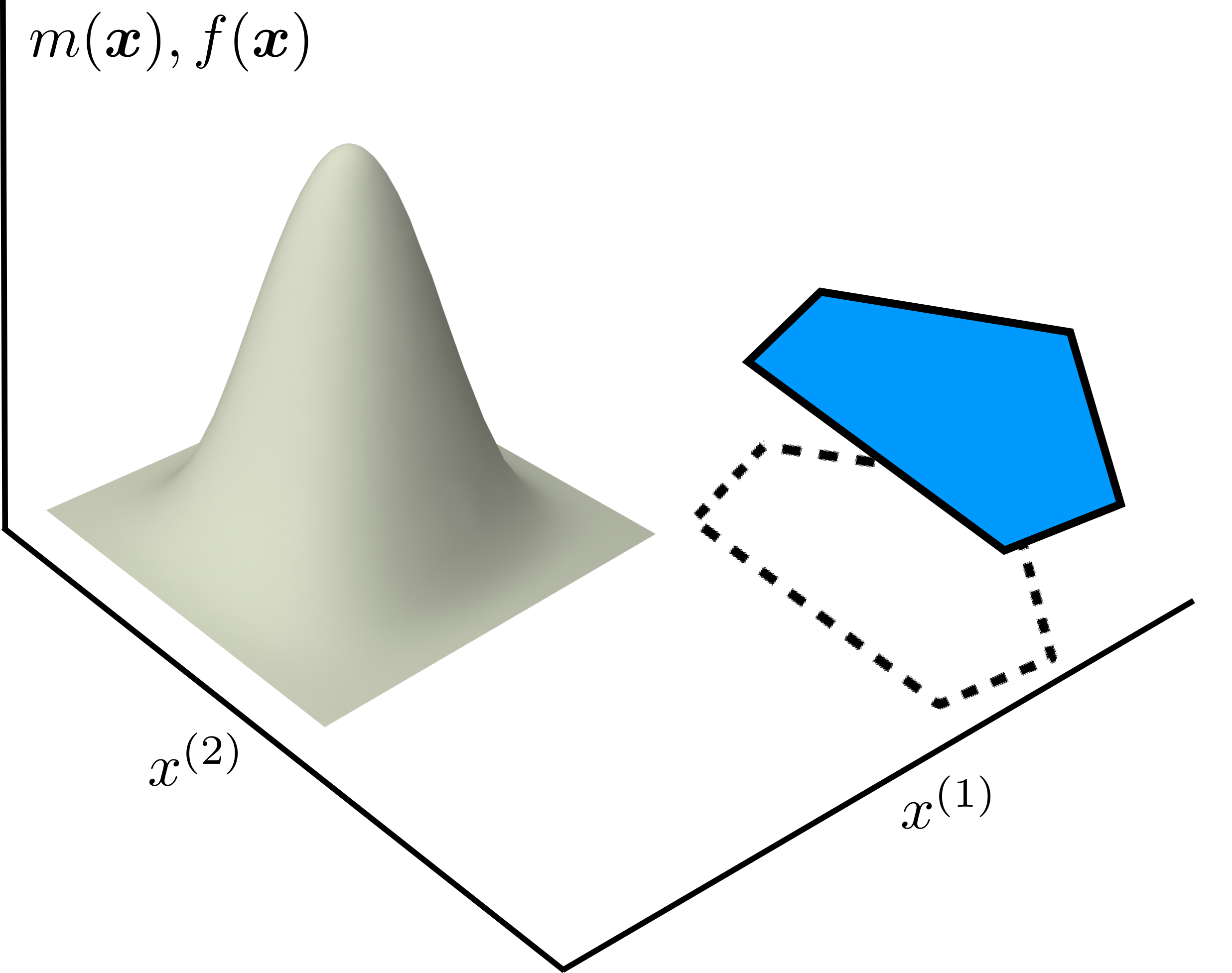} } 
	\hspace{0.1\textwidth}
	\subfloat[][Basis functions of order $r_p = 0$ obtained with the mollifier shown\textbf{} in (a)] {
		\includegraphics[width=0.35 \textwidth]{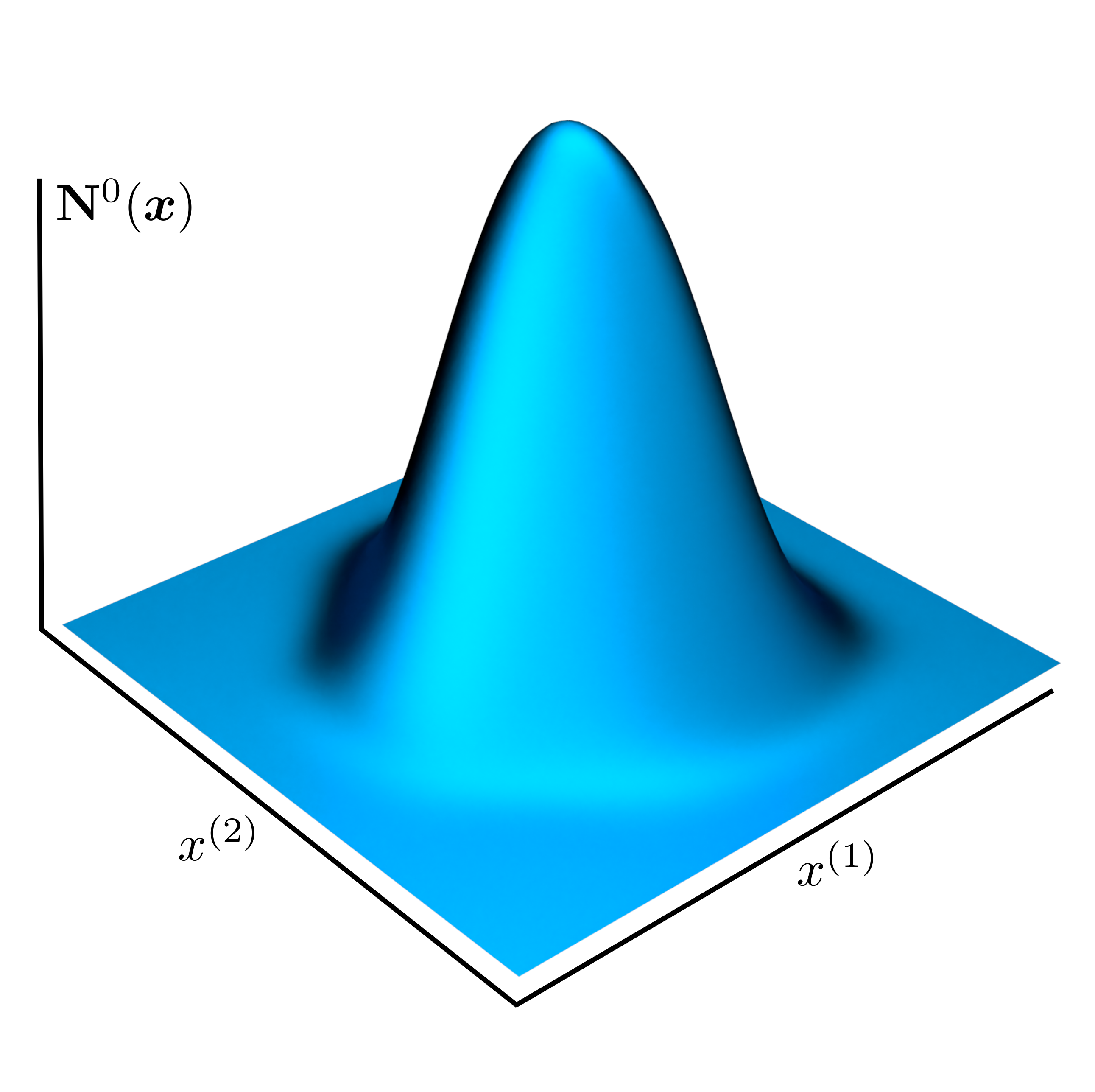} \label{fig:basisSmall}} 
	\\
	\subfloat[][Mollifier and linear local basis functions] {
		\includegraphics[width=0.35 \textwidth]{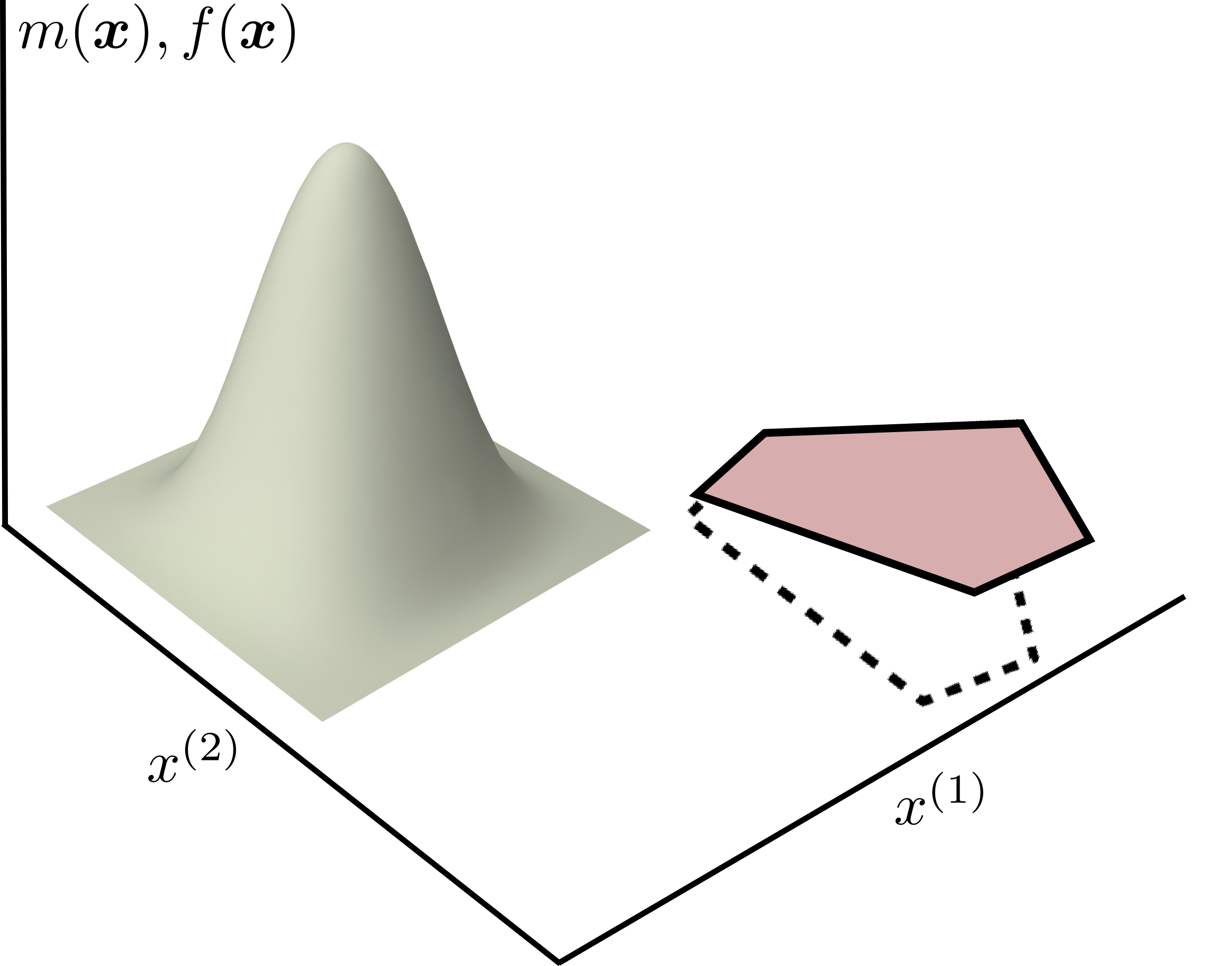} \label{fig:gradSmall}} 
	\hspace{0.1\textwidth}
	\subfloat[][Basis functions of order $r_p = 1$ obtained with the mollifier shown in (c)] {
		\includegraphics[width=0.35 \textwidth]{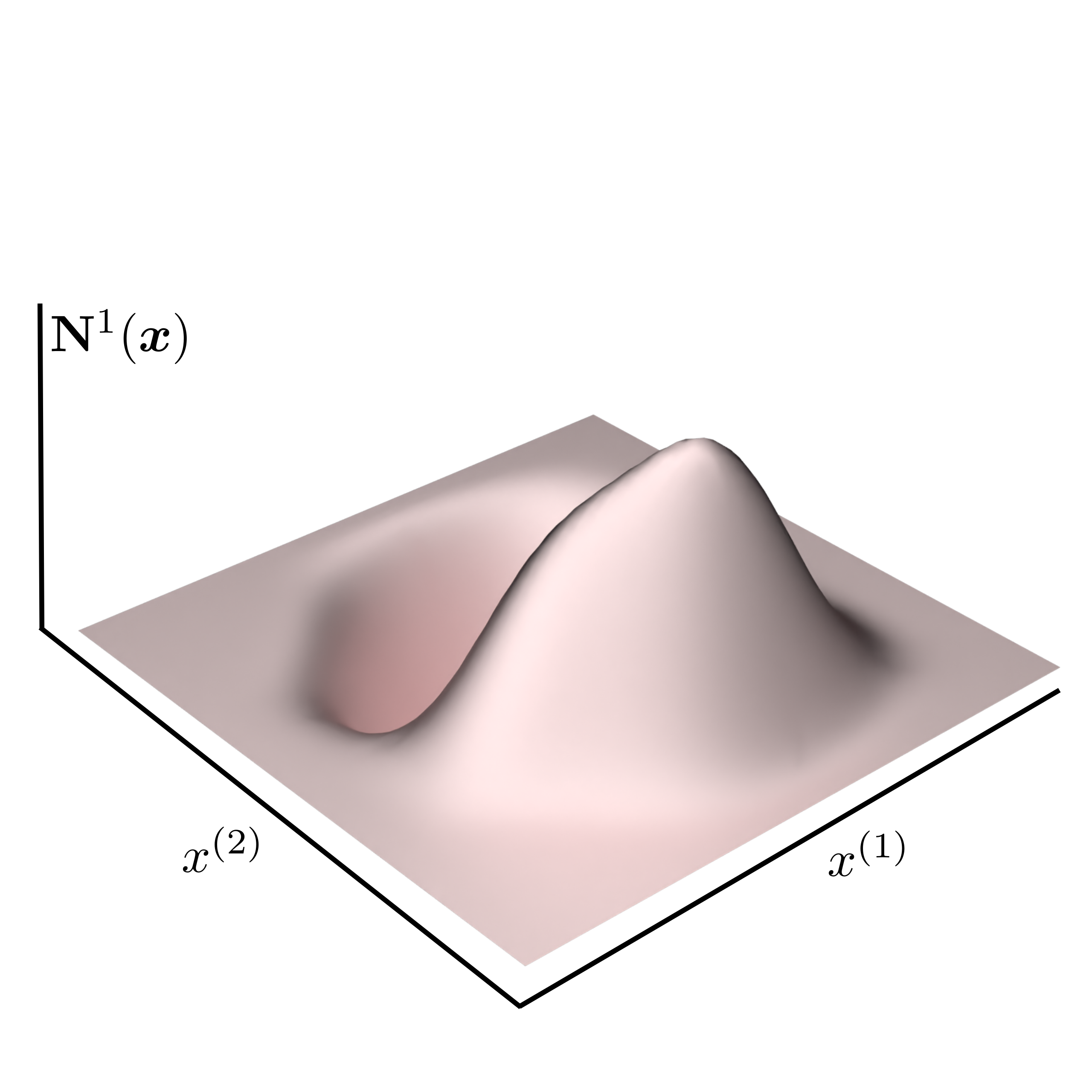} \label{fig:hessSmall}} 
	\caption[Bivariate constant and linear mollified basis functions]{Bivariate constant and linear mollified basis functions with a $C^2$-smooth spline mollifier on a cell~$\omega_i \in \mathbb R^2$. The dashed lines in (a) and (c) indicate the boundary of the cell. \label{fig:basis2D}}
\end{figure}

In each cell $\omega_i$, we consider a set of monomial basis functions centred at the centroid of the cell $\vec x_c$,
\begin{equation}
	\vec p_i(\vec x) = \begin{pmatrix}
		1 & \xi^{(1)} & \xi^{(2)}  & (\xi^{(1)})^2 & (\xi^{(2)})^2 & \dots
	\end{pmatrix}
	\quad \text{  with } \xi^{(j)} =  \frac{2 \left( x^{(j)}- x^{(j)}_c  \right)}{h} \, ,
\end{equation}
where $h$ is the average cell size. In our numerical experiments, $h$ is usually computed as the square root of the average area of the cells in the domain. The bivariate basis functions are obtained by evaluating the convolution integral
\begin{equation} \label{eq:basisMulti}
	\vec N_i(\vec x) =  \int_{\omega_i} m(\vec x - \vec y) \, \vec p_i (\vec y)  \D \vec y  \, .
\end{equation}
As an example, we consider a five-sided polygonal cell $\omega_i$. Mollifying either the constant or linear monomials with a tensor product $C^2$-smooth spline mollifier~\eqref{eq:quartSpline} yields the basis functions presented in Figure~\ref{fig:basis2D}. Likewise, the derivatives of the obtained multivariate basis functions with respect to the coordinate axes can be obtained as in~\eqref{eq:basisDeriv}. 

Next, we focus on the support~$\widehat{\omega} _i$ of the mollified basis functions corresponding to piecewise polynomials of the cell~$\omega_i$. The support~$\widehat{\omega} _i$ is given as the Minkowski sum of the cell~$\omega_i$ with the support of the mollifier~$\Box_{\vec  0} $, that is,
\begin{equation}\label{eq:minkowski}
	\widehat{\omega}_i = \omega_i  \oplus \Box_{\vec  0} = \{ \vec x + \vec y | \,  \vec x \in \omega_i , \vec y \in \Box_{\vec  0} \} \, .
\end{equation}
In our implementation, to obtain the support~$ \widehat{\omega}_i$, we first position the mollifier at the vertices of the cell~$\vec v_{i, j}$, that is,~$\Box_{\vec  v_{i, j}} $. The subscript $j$ is enumerated over the indices of all the vertices of $\omega_i$. We then take the union of the vertices of~$\Box_{\vec  v_{i, j}} $ for all $j$ and combine them using a convex hull algorithm to obtain~$ \widehat{\omega}_i$. For a comprehensive introduction to Minkowski sums, we refer the readers to~\cite{deBerg2010}. 

%--------------------------------------------------------------------------------
\subsection{Basis evaluation at a point}
\label{sec:basis-eval}
%--------------------------------------------------------------------------------
%
In this section, we outline the procedure for evaluating the mollified basis functions at a particular point $\vec x \in \mathbb R^d$. For the univariate case, the convolution integral in~\eqref{eq:basis} can be symbolically computed, resulting in a closed-form expression that can be evaluated at any point $x \in \mathbb R^1$. Such closed-form expressions of the basis functions~\eqref{eq:basis} can be obtained using symbolic programming tools such as Mathematica and SymPy. 

\begin{figure}[h!] 
	\centering
	\subfloat[The mollifier $\Box_{\vec x}$ overlapping some of the cells $\omega_i$]{\includegraphics[width=0.32\textwidth]{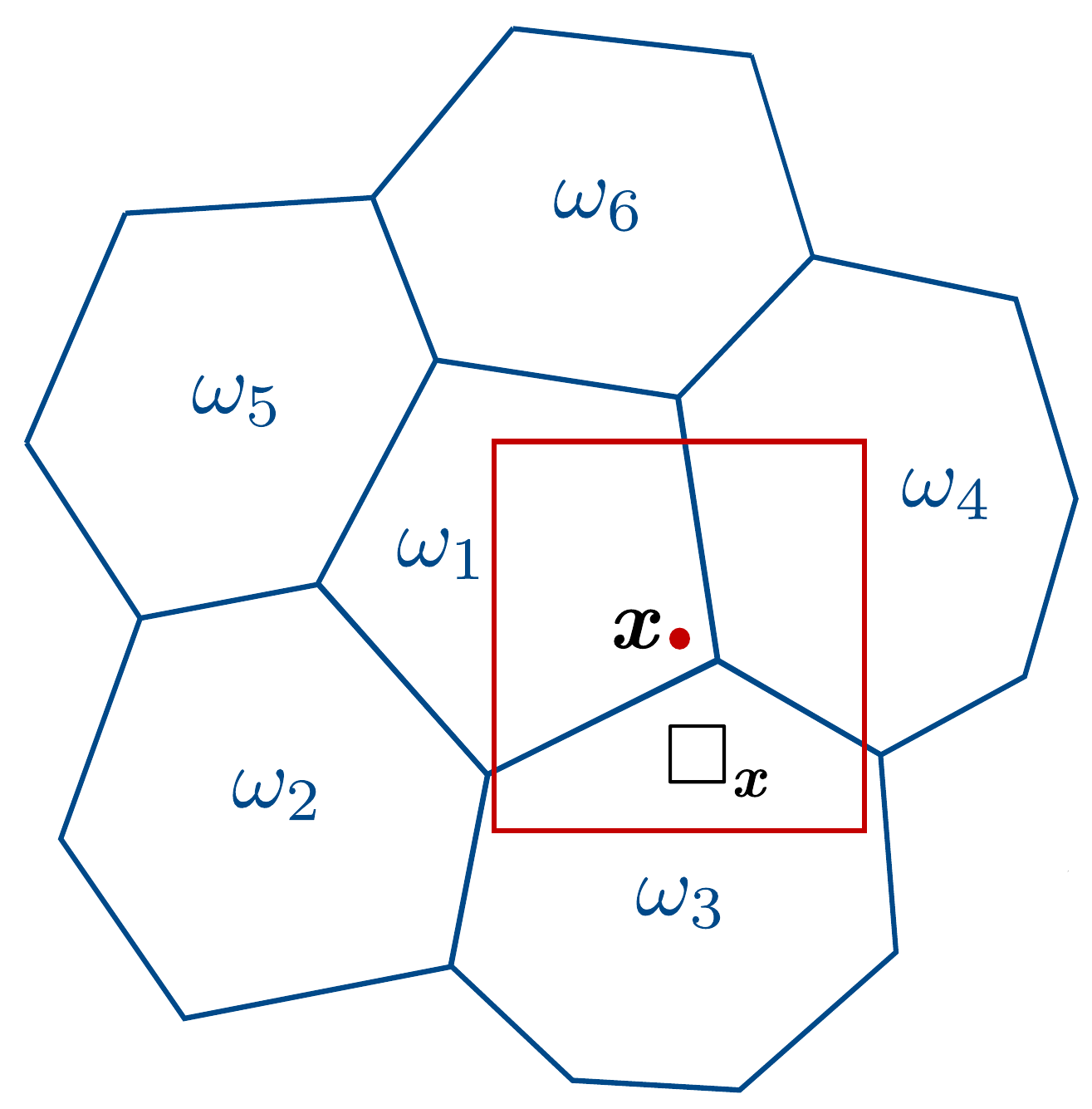} \label{fig:cells}} \qquad 
	\subfloat[$\tau_{i, \vec x}$]{\includegraphics[width=0.28\textwidth]{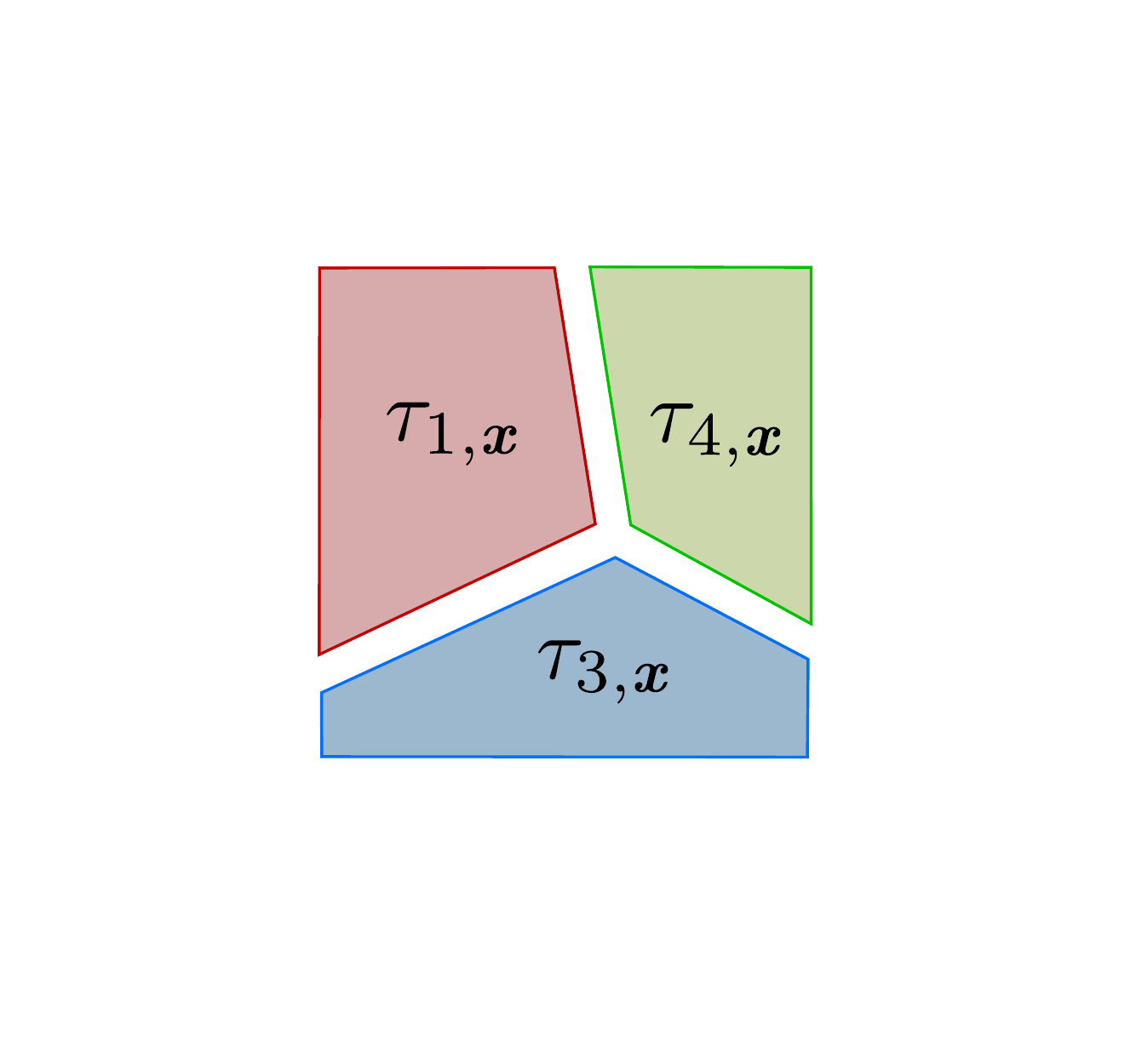} \label{fig:cells-intersect}} \hfill
	\caption[The intersection of the mollifier with the cells $\tau_{i, \vec x}$]{The intersection of the mollifier located at $\vec x$ with the cells $\omega_i$ (a) generates the domain for evaluating the convolution integral $\tau_{i, \vec x}$ (b). }
	\label{fig:intersectionAll}
\end{figure}

When evaluating the multivariate basis functions, the convolution integral in~\eqref{eq:basisMulti} should be numerically computed. To evaluate the basis~$\vec N_i (\vec x)$ at a point~$\vec x \in \mathbb R^2$, we first position the mollifier centred at the evaluation point $\vec x$, defining the support as~$\Box_{\vec x} =   \supp m(\vec x - \vec y)$, as shown in Figure~\ref{fig:intersectionAll}. Both the kernel and polynomial have compact supports, thus, the integrand in~\eqref{eq:basisMulti} is nonzero only within the intersection between the support of the mollifier and the cell. This observation further shrinks the integration domain to
\begin{equation}
	\tau_{i, \vec x}  \coloneqq \Box_{\vec x} \cap  \omega_i  \, ,
\end{equation}
where the intersection~$\tau_i$  is convex because both the cell~$\omega_i$ and the mollifier support~$\Box_{\vec x} $ are convex. The convolution integral used to compute the basis functions is then simplified to
\begin{equation} \label{eq:basisT}
	\vec N_i(\vec x) =  \int_{\tau_{i, \vec x}} m(\vec x - \vec y) \vec p_i (\vec y)  \D \vec y  \, .
\end{equation}

Here, we focus on a two-dimensional example where the local polynomial is defined within a cell~$\omega_i \in \mathbb R^2$ and a polynomial mollifier is employed. To evaluate~\eqref{eq:basisT}, we first triangulate $\tau_{i, \vec x}$ by connecting its edges with its centroid. Gauss quadrature points are then mapped from a reference triangle to each generated triangle to evaluate~\eqref{eq:basisT}. We consider an example using a set of polytopic cells $\left \{ \omega_i \right \}_{i=1}^6$ as shown in Figure~\ref{fig:cells}. To evaluate basis functions at $\vec x$, we place the mollifier support centred at $\vec x$, as shown in Figure~\ref{fig:cells}. The mollifier support $\Box_{\vec x}$ intersects three cells: $\omega_1$, $\omega_3$, and $\omega_4$. The basis functions associated with these cells can be computed by first obtaining their intersection with the mollifier support, namely $\tau_{1, \vec x}$, $\tau_{3, \vec x}$, and $\tau_{4, \vec x}$. Subsequently, we evaluate the convolution integral~\eqref{eq:basisT} over these intersection domains, as shown in Figure~\ref{fig:cells-intersect}.

In the three dimensional case, $\tau_{i, \vec x}$ is a convex polyhedron. The convolution integral over~$\tau_{i, \vec x}$ can be evaluated by tessellating the polyhedron into tetrahedra before distributing the quadrature points. An alternative approach involves the successive application of the divergence theorem to reduce the dimension of integration domains~\cite{mirtich1996fast}. When using a polynomial mollifier, for example, a tensor product spline, the integrand in~\eqref{eq:basisT} is also a polynomial function and can be accurately integrated. For further information on integrating polynomial functions over arbitrary polytopes, interested readers are referred to~\cite{Mousavi2011, chin2020efficient, Antolin2023}. 

%% file: collocation.tex
%--------------------------------------------------------------------------------
\section{Point collocation method}
\label{sec:collocation}
%--------------------------------------------------------------------------------
%
In this section, we outline the discretisation using mollified basis functions for solving PDEs. The smoothness and approximation properties of the mollified basis functions render them suitable for solving PDEs in their strong forms. Our approach assumes that the domain $\Omega$ is discretised into non-overlapping convex polytopes, in which piecewise polynomials are defined. Following the discussion in Section~\ref{sec:basis-eval}, we emphasise that the basis functions and their derivatives can be evaluated at arbitrary points in space, which we later refer to as collocation points. Furthermore, we delve into the considerations involved in determining the collocation points in our numerical studies.
% 
%--------------------------------------------------------------------------------
\subsection{Discretisation}
\label{sec:discretisation}
%--------------------------------------------------------------------------------
%
For simplicity, we consider the Poisson-Dirichlet equation involving a scalar variable $u$ over the domain $\Omega \in \mathbb{R}^d$ with \mbox{$d \in \{ 1, \, 2, \, 3\}$} as a model problem:
\begin{subequations} \label{eq:BVD1D}
	\begin{alignat}{3}
		- \nabla^2 u (\vec x)
		&= s (\vec x)  \qquad && \text{in $\Omega$}   \label{BVDColloc} \, ,
		\\
		u (\vec x)
		&  = \overline{u} \qquad  &&\text{on $\Gamma_D$} \label{BdryColloc} \, ,
	\end{alignat}
\end{subequations}
where $s$ is the source and $\overline{u}$ is the prescribed solution field on the Dirichlet boundary~$\Gamma_D$. The main characteristic of the point collocation method lies in utilising the strong form~\eqref{eq:BVD1D} rather than its weak form. The field variable $u$ is approximated by a linear combination of mollified basis functions, according to the discretisation of the domain $\Omega$ into $n_c$ cells,
\begin{equation} \label{eq:discreteColloc}
	u_h(\vec x) = \sum_{i} \vec N_i( \vec x) \cdot \vec u_i \, .
\end{equation}
In each cell $\omega_i$, the basis functions form a vector consisting of contributions from each monomial. To ensure polynomial reproducibility near the boundary, the domain $\Omega$ is padded with ghost cells with a size proportional to the mollifier width. Substituting $u_h$ as an approximation to $u$ in the Equation \eqref{BVDColloc} yields
\begin{equation} \label{eq:colloc}
	- \sum_{i} \left( \nabla^2  \vec N_i(\vec x) \right) \cdot \vec u_i = s(\vec x) \, .  
\end{equation}
Here, the continuity requirement for the basis functions is at least $C^2$~\cite{Elman2014}. Hence, based on the mollification properties discussed in Section~\ref{sec:mollified}, the minimum continuity for the mollifier $m(\vec x)$ is $C^1$.

We consider a set of points $\{\vec z_j\}_{j=1}^{n_z}$, referred to as collocation points, where Equation \eqref{eq:colloc} is evaluated and the boundary conditions \eqref{BdryColloc} are enforced. We subdivide the point set into two subsets, that is, interior collocation points $\{\vec z^I_k\}_{k=1}^{n_z^I}$ and boundary collocation points $\{\vec z^B_l\}_{l=1}^{n_z^B}$ so that
\begin{equation} \label{eq:collPts}
	\{ \vec z_j \}_{j=1}^{n_z^I + n_z^B} = \left \{\vec z^I_k \right \}_{k=1}^{n_z^I} \cup \left \{\vec z^B_l \right \}_{l=1}^{n_z^B} \, .
\end{equation}
Here, the capital superscripts $I$ and $B$ distinguish points belonging to the interior and boundary set, respectively. Consequently, the total number of collocation points $n_z$ comprises the number of interior points $n_z^I $ and boundary points $n_z^B$, so that
\begin{equation} \label{eq:numColl}
	n_z = n_z^I + n_z^B \, .
\end{equation}
In this approach, we require the total number of collocation points $n_z$ to be greater than, or at least equal to, the total number of the basis functions used in discretising the solution field, denoted as $n_b$. We evaluate the Equation \eqref{eq:colloc} at each interior collocation point $\vec z_k^I$, that is, 
\begin{equation} \label{eq:collocInt}
	- \sum_{i} \left( \nabla^2  \vec N_i (\vec z_k^I ) \right) \cdot \vec u_i = s (\vec z_k^I ) \, ,  
\end{equation}
and enumerate the index $k$ from 1 to the total number of interior points $n_z^I$. The Dirichlet boundary condition is strongly imposed at the boundary collocation points $\vec z_l^B$, that is, 
\begin{equation} \label{eq:collocBdry}
	\sum_{i}  \vec N_i (\vec z_l^B ) \cdot \vec u_i = \overline u \, ,   
\end{equation}
where the index $l$ goes from 1 to the total number of boundary points $n_z^B$. Similarly, for Neumann type of boundary conditions, the gradient of the basis functions $\nabla \vec N_i$ is evaluated at the associated boundary collocation points. 

The two discrete Equations \eqref{eq:collocInt} and \eqref{eq:collocBdry} are compactly expressed in matrix notation as follows
\begin{equation} \label{eq:matColloc}
	\mat C \mat u = \mat s \, .
\end{equation}
Matrix $\mat C$ is a non-square matrix with size $n_z \times n_b$ with $n_z \geq n_b$. Matrix $\mat C$ consists of the two blocks $\mat C^I$ and $\mat C^B$, corresponding to the contributions of the internal and boundary collocation points. Each row of $\mat C^I$ contains the expansion of \eqref{eq:collocInt} when evaluated at an interior point, and $\mat C^B$ consists of the evaluation of \eqref{eq:collocBdry} at boundary collocation points. The entries of $\mat C$ are non-zero only when the collocation point $\vec z_j$ is located within the support of a basis, that is, $\vec z_j \in \supp \vec N_i$. The support of the basis is influenced by the mollifier width, thus, we can deduce that matrix $\mat C$ is denser when the mollifier is wide, and conversely, sparser when the mollifier is narrow. The right-hand size vector $\mat s$ has $n_z$ elements and can be subdivided into $\left(\mat s^I \, \mat u^B \right)^{\trans}$ containing the source terms $s (\vec z_k^I )$ at the interior points and boundary values $\overline u (\vec z_l^B )$ at the boundary points.  

The choice $n_z > n_b$ implies an overdetermined linear system. This system can be solved in a least-square sense by multiplying \eqref{eq:matColloc} with $\mat C^{\top}$, 
\begin{equation} \label{eq:matCollocSq}
	\mat C^{\top} \mat C \mat u = \mat C^{\top} \mat s \quad \Rightarrow  \quad 
	\mat G \mat u = \mat w \, ,
\end{equation}
where $\mat G$ is a square matrix of size $n_b \times n_b$. In the case of $n_z = n_b$, the linear system of Equations \eqref{eq:matColloc} has a solution only when $\mat C$ is non-singular. Therefore, in our implementation we prefer the number of collocation points to be greater than the number of basis functions $n_z > n_b$ to avoid a situation where the system yields no solution. The selection of collocation points in our analysis will be described in the following section. 

% Add a paragraph on the scaling
Our final note concerns the conditioning of the collocation matrix $\mat C$. The mollified basis functions typically involve high-order polynomials due to the order of the local approximants and the mollifier, which potentially leads to poor conditioning of $\mat C$ if untreated. Such conditioning problems become more severe for high-order PDEs. In our implementation, we improve the conditioning of $\mat C$ using two scaling techniques. The first pertains to scaling for the basis functions and the second involves scaling for the derivatives. The first type of scaling is introduced in the description of local approximants~\eqref{eq:scaledMono}, where each monomial is scaled according to the factor $\left(h_c / 2 \right)^{p}$, where $p$ is the monomial degree. This scaling ensures that the maximum value of monomials of any order is the same. The second type of scaling adjusts the $n-$th derivative of the basis function using factor $\left( h_m \right)^{n}$ to account for magnitude discrepancies across the derivatives. Evidently, the inverse of the second scaling factor should be applied upon acquiring the numerical solution of the linear system. 

%--------------------------------------------------------------------------------
\subsection{Spatial distribution of collocation points}
\label{sec:collocation-select}
%--------------------------------------------------------------------------------
%
The entries in the linear system \eqref{eq:matColloc} depend on the location of the collocation points $\vec z_j$ where the basis functions and their second derivatives are evaluated. Various strategies have been proposed for selecting the position of the collocation points. When piecewise polynomial approximants are used, Gauss quadrature can be devised as collocation points~\cite{deBoor1973, Finlayson2013}. For mesh-free approximation schemes, such as the reproducing kernel particle method (RKPM) and local max-ent approximation, collocation points can be chosen either similarly to the points representing the basis functions~\cite{Aluru2000, greco2020}, or differently~\cite{hu2007, Fan2018}. In the context of isogeometric collocation using B-spline and NURBS basis functions, Greville and Demko abscissae are employed as the locations for collocating the strong form of the governing Equation \eqref{eq:colloc}~\cite{Johnson2005, Auricchio2010, Manni2015, Wang2021}, which exploits the knot structure of B-splines. Moreover, recent studies have reported the existence of collocation points with improved convergence properties similar to the Galerkin schemes~\cite{Anitescu2015, Gomez2016, montardini2017}. 

As mentioned in Section~\ref{sec:mollified-basis}, mollified basis functions are obtained by convolving piecewise polynomials defined over a polytopic cell with a mollifier. The support of the basis functions generally do not coincide with the cell boundaries. Therefore, it is challenging to exploit a specific structure to find optimal coordinates for collocation points when discretising with mollified basis functions. This is unlike the case with B-splines, where the tensor-product structure can be used to identify optimal or privileged collocation points. Nonetheless, this work does not pursue the identification of such optimal coordinates for collocation points. Furthermore, collocating at Voronoi seeds, similar to the nodes in meshfree methods, leads to underdetermined  system matrices because each cell has multiple associated basis functions. Therefore, we explore several strategies for distributing the collocation points, including uniform, Gauss quadrature, and quasi-random point distributions. 

To ensure stability, we choose the number of collocation points so that the linear system matrix is overdetermined. This allows for a less stringent approach to selecting the position of collocation points. In our approach, the number of collocation points, including those on the boundary, can be described by expanding the interior points $n_z^I$ in terms of the number of cells 
\begin{align} \label{eq:num-colloc-general}
	n_z = n_z^I + n_z^B = \beta \, n_c + n_z^B \, \geq n_b  \, .
\end{align} 
Here, the factor $\beta \geq 1$ accounts for the total number of monomial basis functions in each cell, which directly corresponds to the polynomial order $r_p$. The number of boundary collocation points $n_z^B$ is explicitly separated in~\eqref{eq:num-colloc-general} so that $\beta$ can conveniently be chosen as an integer. For example, in the set of basis functions shown in~Figure~\ref{fig:basisAllUnif}, the total number of basis functions is $(n_c + 2) \, |\vec p|$, where $n_c = 6$ and $| \vec p | = r_p + 1 = 4$. Moreover, the addition with two accounts for the appended ghost cells. Therefore, the total number of basis functions involved is $n_b = 32$. In this case, the closest integer $\beta$ is selected so that~\eqref{eq:num-colloc-general} is satisfied, which is $\beta = 6$. 

In this paper, we consider three methods for distributing collocation points: Gauss quadrature, uniform, and quasi-random point distribution. In the univariate setting, obtaining these point distributions is straightforward, as depicted in Figure~\ref{fig:coll-point-dist}. In the multivariate setting, we can obtain the Gauss point arrangement by first tessellating the convex polytopic cells into simplices (that is, triangles or tetrahedra in two or three dimensions, respectively). Gauss points of order $\gamma$ are then mapped from the reference simplex. The uniform collocation point distribution can be obtained by placing the points in an equidistant manner in each coordinate axis and  subsequently taking their tensor products in higher dimensions. Note that we only consider uniformly distributed collocation points for numerical examples with tensor product domains. Lastly, the quasi-random point collocation distribution is obtained by introducing perturbation to an initial point distribution, which is usually chosen as uniform. The random perturbation $\epsilon \sim \set U(-\sigma, \, \sigma)$ is generated in each coordinate axis from a univariate uniform density function with maximum perturbation $\sigma$. The coordinates of collocation points are detailed for each example in Section~\ref{sec:example}.

%% file: examples.tex
%--------------------------------------------------------------------------------
\section{Examples}
\label{sec:example}
%--------------------------------------------------------------------------------
%

In this section, we present several numerical experiments with increasing complexity to investigate the convergence property of the proposed mollified-collocation approach. We initially investigate the performance of the mollified-collocation scheme in solving one-dimensional Poisson and biharmonic problems. The effects of the piecewise polynomial degree, mollifier, and collocation points on the convergence are studied. We then proceed with the numerical study of a two-dimensional elastic plate, a plate with a hole, a plate bending, and a three-dimensional cube. In all the examples, convergence is studied using the relative discretisation errors as in~\cite{greco2020}. The relative error between field $\vec v$ and its numerical approximation $\vec v_h$ can be defined as
\begin{equation}
	\label{key}	e(\vec v, \vec v_h) = \left( \frac{\sum_k^{n_z} (\vec v^k- \vec v_h^k) (\vec v^k - \vec v_h^k)^{\trans}}{\sum_k^{n_z} \vec v^k (\vec v^k)^{\trans}}  \right)^{\frac{1}{2}} \, ,
\end{equation}
where $\vec v^k$ and $\vec v_h^k$ are the evaluation of the field and its approximation at collocation point $\vec z_k$.  

%--------------------------------------------------------------------------------
\subsection{One-dimensional examples} 
\label{sec:onedim}
%--------------------------------------------------------------------------------

%--------------------------------------------------------------------------------
\subsubsection{One-dimensional Poisson problem} 
\label{sec:onedimPoi}
%--------------------------------------------------------------------------------

%
As a first example we consider the solution of the one-dimensional Poisson-Dirichlet problem $- \D^2 u(x) / \D x^2  =  s(x)$ on the domain $\Omega = (0, \,1 ) \in \mathbb R^1$. The source term~$s(x)$ is chosen such that the solution is equal to $u(x) = \sin (3 \pi x)$. To construct the mollified basis functions, we first consider the decomposition of the domain $\Omega$ into a set of non-overlapping cells $\{\omega_i\}_{i=1}^{n_c}$ of size $h_{c,\, i}$ where the piecewise polynomial basis functions $\vec p_i (x)$ are defined. We initially consider the coarse $n_c = 6$ non-uniform cells with size of $h_{c,1} = 0.15$, $h_{c,2} = 0.2$, $h_{c,3} = 0.15$, $h_{c,4} = 0.15$, $h_{c,5} = 0.2$ and $h_{c,6} = 0.15$. We then consider refinement through the bisectioning of each cell. The convolution integral is evaluated by assuming a mollifier function $m(x)$ with support size of 
\begin{equation}
	\label{eq:moll-width}
	h_m = 2 \left(\max_j  h_{c,j} \right) \, .
\end{equation} 
The closed form of the basis functions can be obtained by analytically evaluating the convolution integral involving the piecewise polynomials and the mollifier. To ensure completeness at the boundary, one ghost cell layer of size $h_m$ is padded to each end of the domain. In the following set of experiments, we study the influence of collocation point distribution, local polynomial basis $\vec p_i(x)$, and the mollifier~$m(x)$ on the convergence of the solution approximation. 

%Effect of collocation points distribution
We first study the influence of three types of point distribution: uniform, Gauss quadrature, and quasi-random, on the error convergence. Even though the knot-based abscissae are applicable in the one-dimensional case, they are not extendable for higher dimensions with unstructured polytopic partitions, and therefore, are not considered in the present analysis. The uniform point distribution is obtained by placing the collocation points over the domain $\Omega$ in an equidistant manner without considering the cell boundaries. By contrast, the Gauss quadrature points are distributed by mapping the standard Gauss quadrature from the parametric domain onto each cell. Finally, the quasi-random point distribution is obtained by adding a small random perturbation sampled from a uniform distribution $\epsilon \sim \set U(-\sigma, \, \sigma)$ to the uniformly distributed interior collocation points. Here the parameter $\sigma$ is chosen as 10\% of the equidistant spacing of the uniform configuration. 

\begin{figure}[]
	\centering
	\includegraphics[scale=0.35]{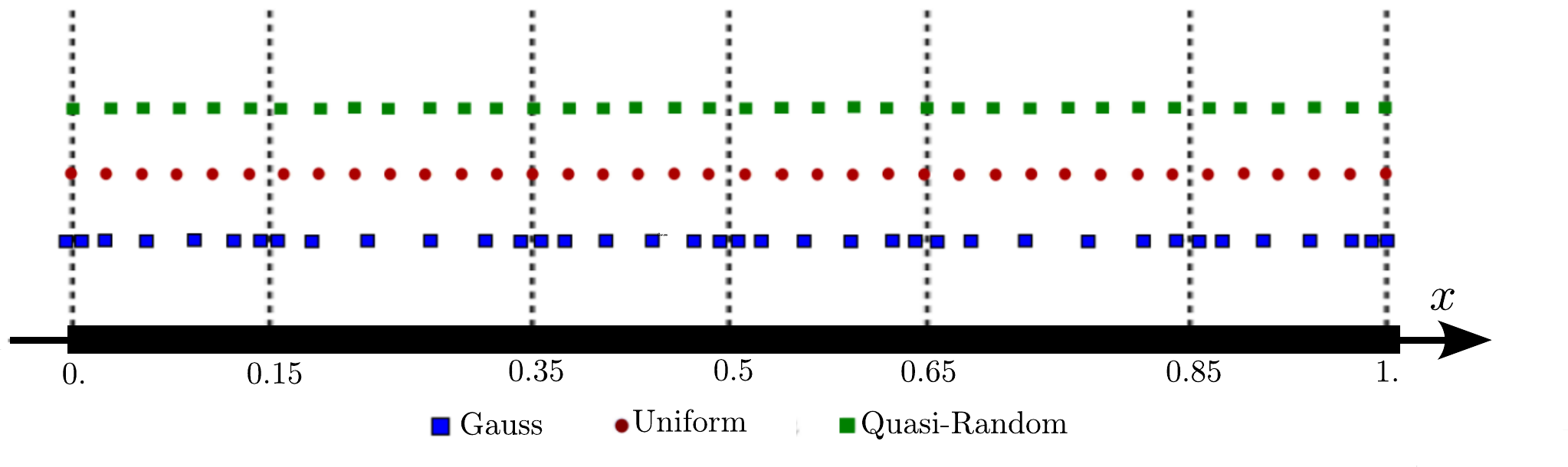}
	\caption{One-dimensional collocation point distributions.}
	\label{fig:coll-point-dist}
\end{figure}

The total number of collocation points $n_z$ should not be less than the number of basis functions $n_b$ to ensure that the system matrix is not underdetermined. In our one dimensional example, the total number of basis functions involved in the computation $n_b$ depends on the number of cells $n_c$ and the order of the local polynomial $r_p$
\begin{align} \label{eq:coeff-num-1d} 
	n_b =  \left( n_c + 2 \right) \, (r_p + 1) \, .
\end{align} 
The above expression considers the contribution from one layer of ghost cell on each side. We consider the total number of collocation points $n_z$ following~\eqref{eq:numColl} 
\begin{align} \label{eq:coll-num} 
	n_z = \left(\beta \, n_c \right) + 2 \, ,
\end{align} 
where the term inside the bracket resembles the internal collocation points. For consistency, we use $\beta = 6$ in the examples throughout this section following the minimum number of unknown coefficients imposed by the cubic polynomial basis $r_p = 3$ at the coarsest level $n_c = 6$. The constant two in \eqref{eq:coll-num} accounts for the boundary collocation points where the Dirichlet boundary conditions are imposed. Figure~\ref{fig:coll-point-dist} depicts the three distributions of the collocation points for the one-dimensional test case with cubic polynomials $r_p = 3$ and $n_c = 6$ cells. 

\begin{figure}[h!]
	\centering
	\subfloat[][$L^2$-norm error \label{fig:1d coll L2}]{
		\includegraphics[scale=0.5]{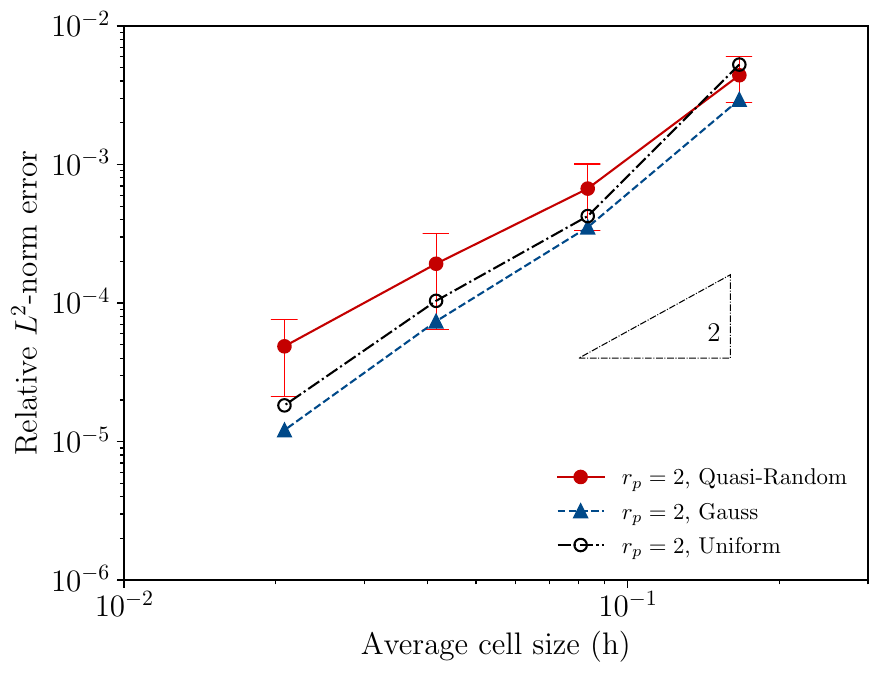}
	}
	\hspace{0.05\textwidth}
	\subfloat[][$H^1$-seminorm error \label{fig:1d coll H1}]{
		\includegraphics[scale=0.5]{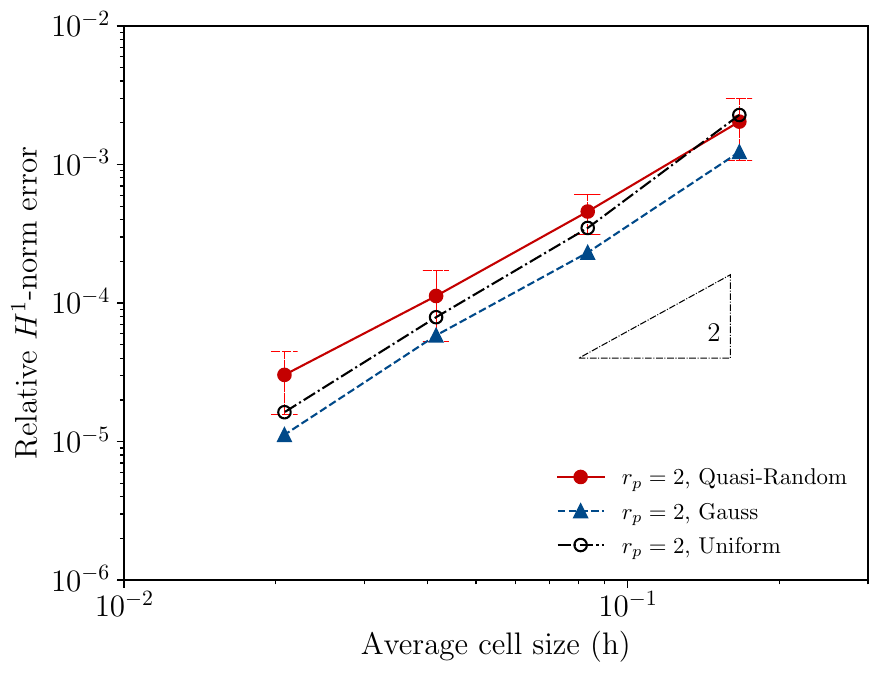}
	}
	\caption[]{One-dimensional Poisson problem. Convergence with normalised quadratic B-spline mollifier and quadratic polynomial basis for three collocation point distributions: Uniform, Gauss quadrature, and quasi-random.}
	\label{fig:coll-dist-conv}
\end{figure}

Figure~\ref{fig:coll-dist-conv} shows the error convergence in both the $L^2$-norm and $H^1$-seminorm when using a normalised quadratic B-spline mollifier and quadratic polynomial order $r_p = 2$. In particular, for the perturbed point distribution, the random perturbation $\epsilon$ is first obtained for each interior collocation point, generating a set of quasi-random points. We obtain such a point distribution one hundred times, from which the mean and standard deviation of solution errors can be obtained. Figure~\ref{fig:coll-dist-conv} shows that the three collocation point distributions yield a convergence rate of $2$ in both the $L^2$-norm and $H^1$-seminorm, which aligns with the collocation methods using IGA and max-ent basis functions~\cite{Auricchio2010, greco2020}. In Figure~\ref{fig:coll-dist-conv}, the three point distributions yield only small differences in the convergence constants. Furthermore, the quasi-random point arrangement has the highest mean error among the three point distributions considered. On the other hand, the quasi-random point distribution has the lowest convergence rate of mean errors among the three point distributions. Moreover, for this distribution, the standard deviation of the errors increases as refinement progresses. 
%to follow the cell shape and the breakpoints of the basis functions.

\begin{figure}[h!]
	\centering
	\subfloat[][$L^2$-norm error \label{fig:1d-poly-L2}]{
		\includegraphics[scale=0.5]{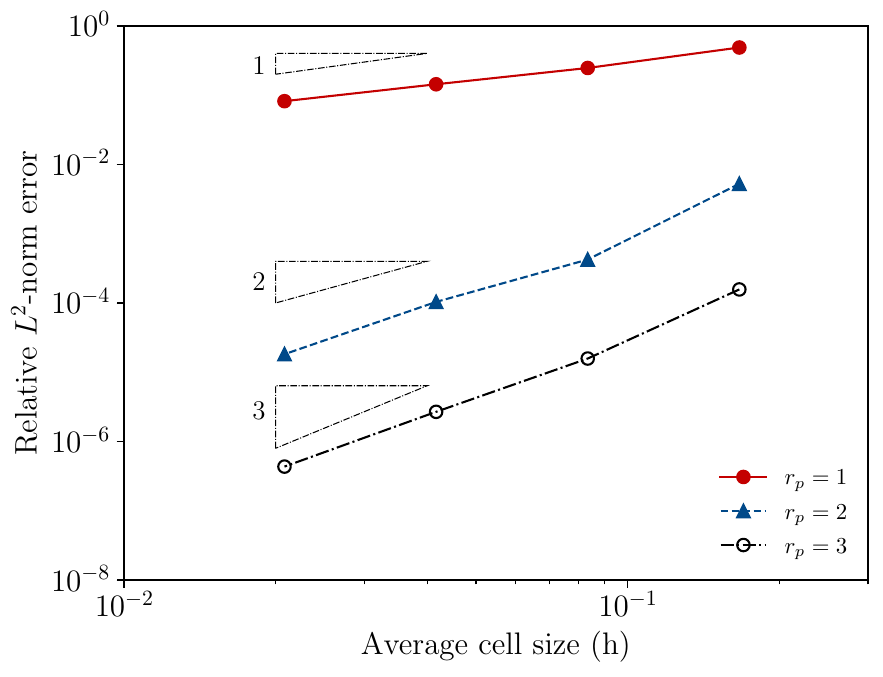}
	}
	\hspace{0.05\textwidth}
	\subfloat[][$H^1$-seminorm error\label{fig:1d-poly-H1}]{
		\includegraphics[scale=0.5]{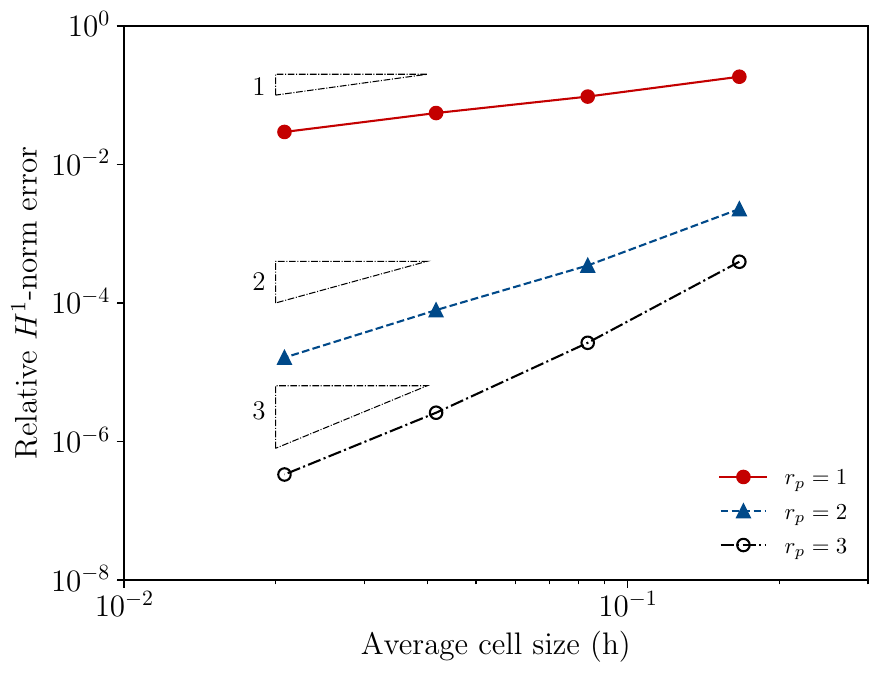}
	}
	\caption[]{One-dimensional Poisson problem. Convergence with normalised quadratic B-spline mollifier and piecewise polynomial of order $r_p = \{1, 2, 3\}$.}
	\label{fig:poly-study}
\end{figure}

%Effect of polynomial order
The second experiment aims to study the influence of the local polynomial order~$r_p \in \{ 1,2,3\}$ on the error convergence of the numerical approximation. We consider a normalised quadratic B-spline mollifier with width chosen according to~\eqref{eq:moll-width}. In this example, we choose uniformly distributed collocation points with factor $\beta = 6$ and spacing $1 / (n_z + 1)$. Figure~\ref{fig:poly-study} presents the error convergence for each polynomial degrees $r_p$ in both the $L^2$-norm and $H^1$-seminorm. The convergence rates in the $L^2$ norm are $0.80$, $2.27$, and $2.58$ for polynomial orders $r_p = 1$, $r_p = 2$, and $r_p = 3$, respectively. A similar trend can be observed in the $H^1$-seminorm, where the average convergence rates are $0.85$, $2.21$, and $3.16$ for polynomial orders $r_p = 1$, $r_p = 2$, and $r_p = 3$, respectively. These results suggest that the proposed approach converges with rate $r_p$ in both error norms with lower rates for the odd polynomial order, particularly in the $L^2$-norm, which resembles the findings reported in previous studies~\cite{Auricchio2010, greco2020}.

\begin{figure}[h!]
	\centering
	\subfloat[][$L^2$-norm error \label{fig:1d-moll-width-L2}]{
		\includegraphics[scale=0.5]{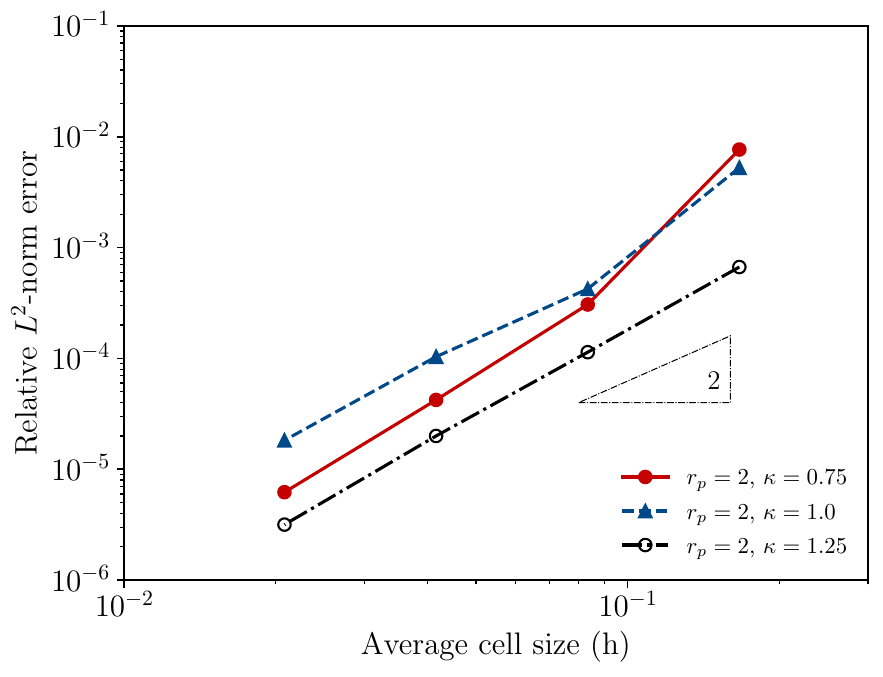}
	}
	\hspace{0.05\textwidth}
	\subfloat[][$H^1$-seminorm error \label{fig:1d-moll-width-H1}]{
		\includegraphics[scale=0.5]{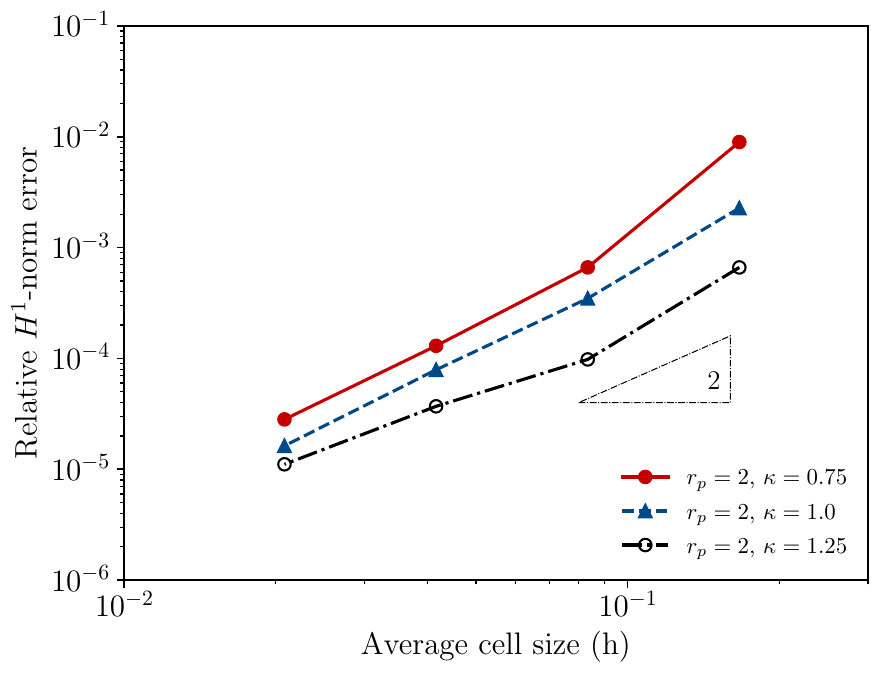}
	}
	\caption[]{One-dimensional Poisson problem. Convergence with normalised quadratic B-spline mollifier of width factor $\kappa = \{ 0.75, 1., 1.25\}$ and quadratic piecewise polynomial. }
	\label{fig:moll-width-study}
\end{figure}

%Effect of mollifier
Subsequently, we investigate the influence of mollifier width and smoothness on the error convergence. We use in this example uniformly distributed collocation points with factor $\beta = 6$ and a piecewise polynomial of order $r_p = 2$. For the first case, we modify the width of a quadratic B-spline mollifier through scaling factor $\kappa$, that is,  
\begin{equation}
	h_m =  2  \kappa \left ( \max_j h_{c,j} \right )  \quad \text{with } \kappa \in \{0.75, \,1., \,1.25 \} \, .  
\end{equation} 
The increase in mollifier width leads to larger support of the basis functions. This implies that more basis functions have non-zero values at a collocation point which leads to a denser matrix system. Figure~\ref{fig:moll-width-study} shows that changing the mollifier width leads to only a slight difference in the error magnitude while keeping the convergence rate as $r_p = 2$ for both the $L^2$-norm and $H^1$-seminorm. 

\begin{figure}[h!]
	\centering
	\subfloat[][$L^2$-norm error \label{fig:1d-moll-order-L2}]{
		\includegraphics[scale=0.5]{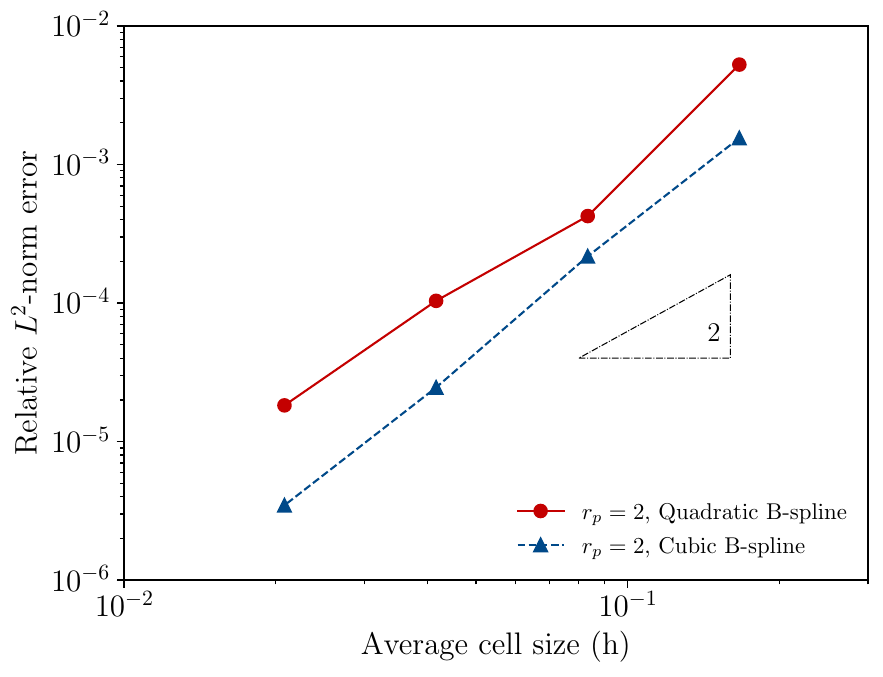}
	}
	\hspace{0.05\textwidth}
	\subfloat[][$H^1$-seminorm error \label{fig:1d-moll-order-H1}]{
		\includegraphics[scale=0.5]{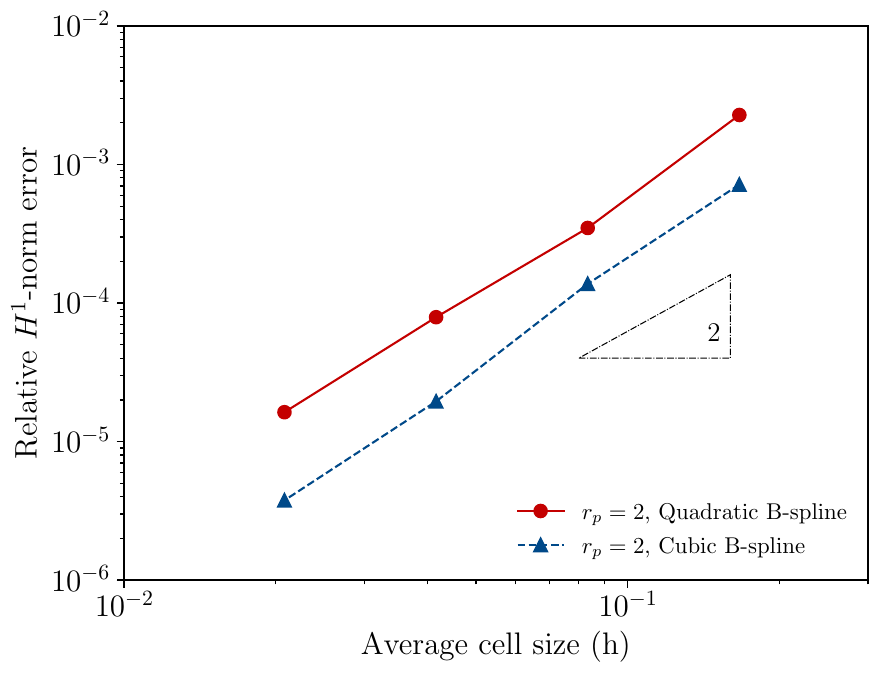}
	}
	\caption[]{One-dimensional Poisson problem. Convergence with normalised quadratic and cubic B-spline mollifier and quadratic piecewise polynomial $(r_p = 2)$. }
	\label{fig:moll-order-study}
\end{figure}

We next investigate the effect of mollifier smoothness on the convergence of the approximation error. Here, we consider quadratic and cubic B-spline mollifiers which are~$C^{1}$- and~$C^{2}$-smooth, respectively. The obtained mollified basis functions are~$C^{2}$- and $C^{3}$-smooth, respectively. This aligns with the smoothness requirement from the PDE that the solution has to be at least $C^2$-continuous. The mollifier width is chosen according to~\eqref{eq:moll-width}, i.e, $\kappa = 1$. Figure~\ref{fig:moll-order-study} displays the convergence plot for both the $L^2$-norm and $H^1$-seminorm errors. It is evident that higher mollifier smoothness improves the convergence constants while keeping a convergence rate of $r_p = 2$. 

Finally, we address the cost comparison between the proposed mollified-collocation method and the mollified-FEM~\cite{Febrianto2021} by estimating the number of point evaluations required to construct the linear system. We consider a quadratic polynomial with $r_p = 2$ over 48 cells, mollified with a quadratic B-spline mollifier. The resulting basis functions have, as described in~\eqref{eq:basis}, monomials up to order 5. Consequently, the stiffness integrand in the mollified-FEM has monomials up to order 8, necessitating 5 Gauss quadrature points on each non-overlapping polynomial segment, thus yielding a total of 240 point evaluations. In contrast, the proposed mollified-collocation method only requires evaluation at a minimum 150 collocation points. This discrepancy is further accentuated in higher dimensions and when using a finer discretisation. The number of quadrature points in mollified-FEM can be somewhat reduced by using the variationally consistent integration (VCI) approach. VCI requires the solution of a small dense matrix to derive the corrective coefficients in each cell.
% estimation of total point evaluation (collocation vs FEM)
% increasing polynomial order
%1  100 vs 192 100%
%2 150 vs 240  60%
%3 200 vs 288 44%
% finer discretisation
%12 42 vs 60 43%
%24 78 vs 120 54%
%48 150 vs 240 60%

%--------------------------------------------------------------------------------
\subsubsection{One-dimensional biharmonic problem} 
\label{sec:oneDbiharmonic}
%--------------------------------------------------------------------------------

We next consider the one-dimensional biharmonic problem $\D^4 u(x) / \D x^4  = s(x)$ over the domain $\Omega = (0, \,1 ) \in \mathbb R^1$. The source term $s(x)$ is chosen such that the solution is $u(x) = \sin (3 \pi x)$. Both the value $u$ and the first derivative $u'$ are prescribed at the boundaries $x=0$ and $x = 1$. In this example, the domain $\Omega$ is uniformly discretised into $n_c$ cells of size $h_c = 1 / n_c$. In each cell, we consider a polynomial of order $r_p \in \{ 5, \, 6\}$ consisting of $r_p + 1$ monomial basis and their respective coefficients. Because of the higher derivatives that appear in the PDE, the smoothness requirement for the solution $u(x)$ is $C^4$, which requires that the mollifier should be at least $C^3$ smooth. Hence, we consider a $C^3$ spline mollifier 
\begin{equation} \label{eq:octicSpline}
	m(x) =  
	\begin{cases} 
		\frac{315}{128 h_m} \left(  1 - 16 \left( \frac{x}{h_m} \right)^2 + 96 \left( \frac{x}{h_m} \right)^4 - 256 \left( \frac{x}{h_m} \right)^6 + 256 \left( \frac{x}{h_m} \right)^8  \right) & \text{  if  } |x| < \frac{h_m}{2}\\ 
		0 & \text{  if  } |x| \geq \frac{h_m}{2}\end{cases}
\end{equation}
where the mollifier width $h_m$ is chosen to be twice the cell size $h_m = 2 h_c$. We uniformly distribute $n_z$ collocation points with the factor $\beta$ according to~\eqref{eq:coll-num}. We analyse the effect of the  total number of collocation points in terms of the factor $\beta$, where $\beta = \{8, 10\}$ are compared for the fifth-order polynomial and $\beta = \{10, 12\}$ are considered for the sixth-order polynomial. These factors are larger than the ones specified in the Poisson examples (Section ~\ref{sec:onedimPoi}) because of the higher order polynomials involved. 

\begin{figure}[]
	\centering
	\subfloat[][$L^2$-norm error \label{fig:biharmonic-L2}]{
		\includegraphics[scale=0.5]{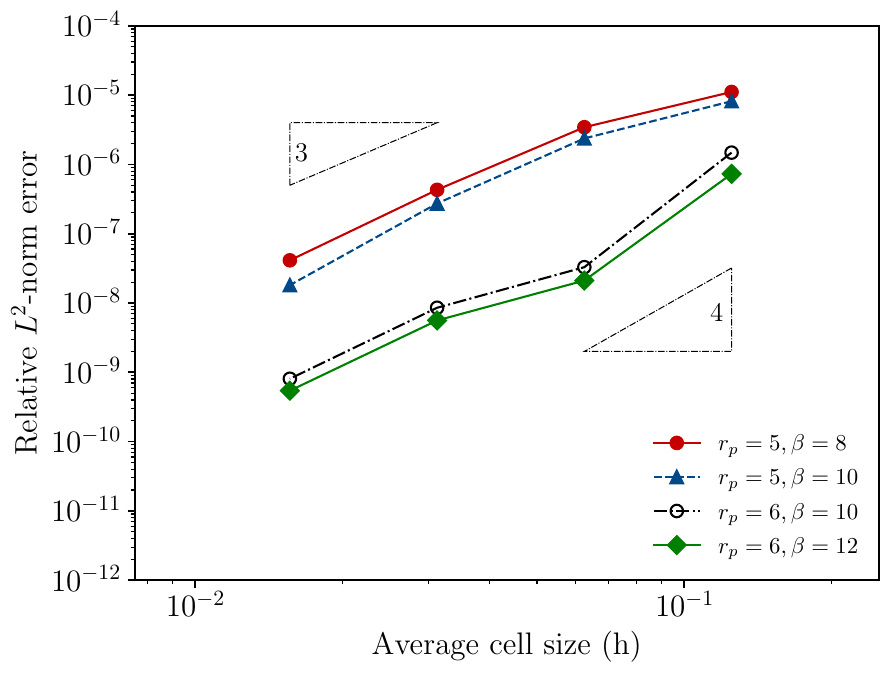}
	}
	\hspace{0.05\textwidth}
	\subfloat[][$H^1$-seminorm error \label{fig:biharmonic-H1}]{
		\includegraphics[scale=0.5]{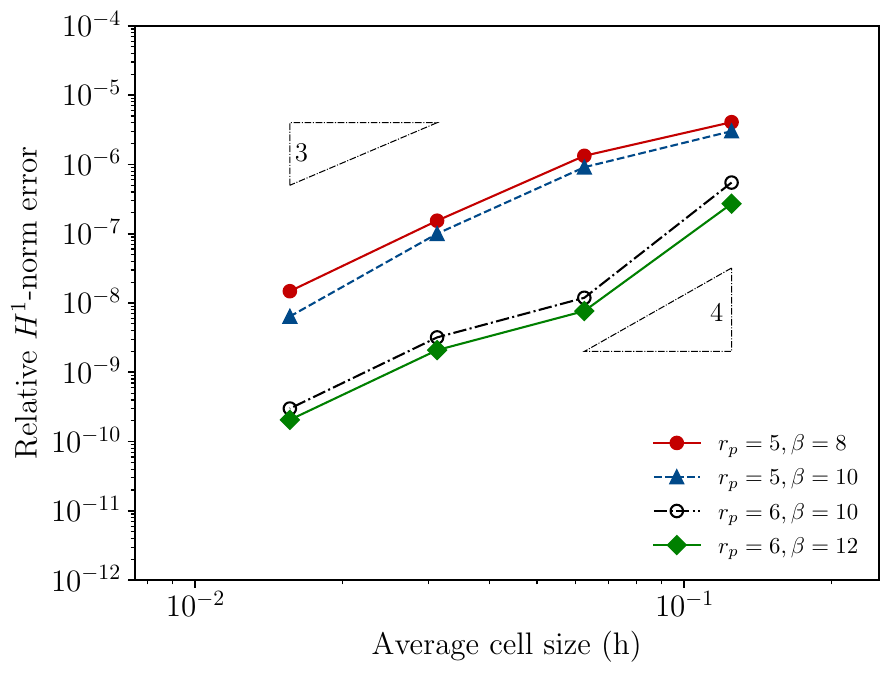}
	}
	\caption[]{One-dimensional biharmonic problem. Convergence with the $C^3$ spline mollifier with a piecewise polynomial of order $r_p = \{5, 6\}$. }
	\label{fig:biharmonic}
\end{figure}

The convergence of relative errors in the $L^2$-norm and $H^1$-seminorm are shown in Figure~\ref{fig:biharmonic}. In the $L^2$-norm, the error converges approximately with rate $r_p - 2$ for both $r_p = 5$ and $r_p = 6$. In the $H^1$-seminorm, a convergence rate of $r_p - 2$ can also be observed for both polynomial orders. The number of collocation points, dictated by factor $\beta$, affects the convergence constants in both norms. For the same polynomial order, a higher $\beta$ indicates more collocation points, which yields lower convergence constants. Morever, a higher number of collocation points improves the conditioning of the system matrix~\eqref{eq:matColloc}.  

%
%--------------------------------------------------------------------------------          
\subsection{Two-dimensional examples}
\label{sec:twodim}
%--------------------------------------------------------------------------------

%--------------------------------------------------------------------------------          
\subsubsection{Two-dimensional elasticity}
\label{sec:plate}
\begin{figure}[]
	\centering
	\subfloat[][$n_c = 16$ \label{fig:cell16U}]{
		\includegraphics[scale=0.45]{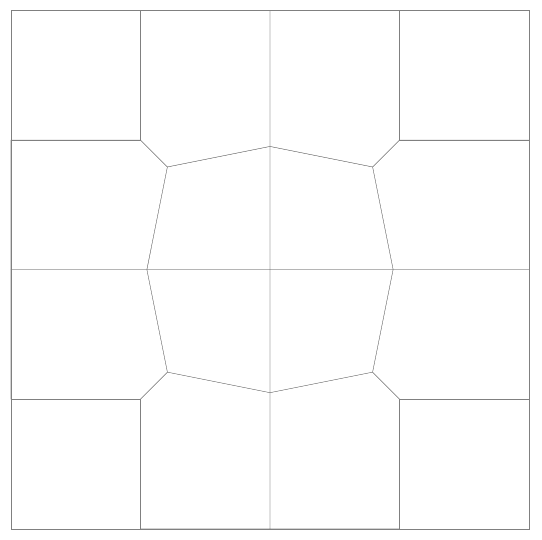}
	}
	\hspace{0.03\textwidth}
	\subfloat[][$n_c = 64$ \label{fig:cell64U}]{
		\includegraphics[scale=0.45]{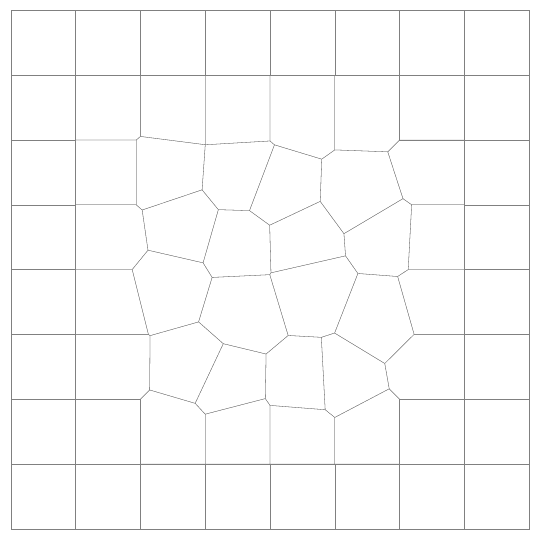}
	}
	\hspace{0.03\textwidth}
	\subfloat[][$n_c = 256$ \label{fig:cell256U}]{
		\includegraphics[scale=0.45]{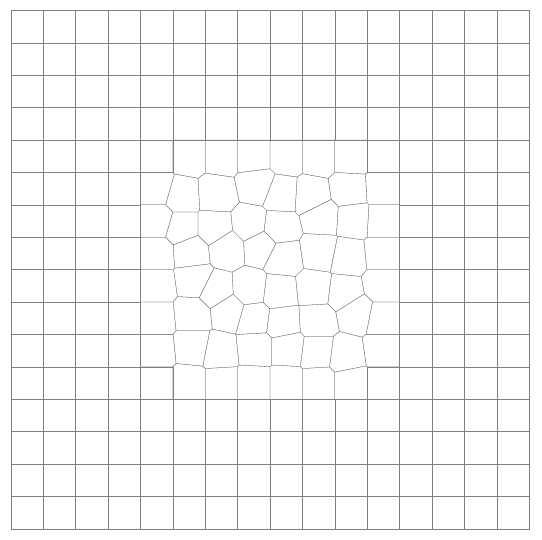}
	}
	\caption[]{Voronoi diagram consisting of $n_c = 16$ (a), $n_c=64$ (b), and $n_c = 256$ (c). The Voronoi seeds are obtained from~\cite{Febrianto2021}.}
	\label{fig:cells2d}
\end{figure}

%--------------------------------------------------------------------------------
%
We next consider the linear elasticity problem $\nabla \cdot \vec \sigma (\vec x) = \vec b (\vec x)$ on the two-dimensional square domain $\Omega = (0,1) \times (0,1)$. The material of the plate has Young’s modulus $E = 1000$ and its Poisson’s ratio is $\nu = 0.3$. The body force $\vec b(\vec x)$ is chosen such that the solution equals to 
\begin{equation}
	u^{(1)} = u^{(2)} = \sin(\pi x^{(1)}) \sin(\pi x^{(2)}) \, . 
\end{equation}
The domain $\Omega$ is partitioned into $n_c$ cells using the Voronoi diagram of $n_c$ non-uniformly distributed seeds. The cells corresponding to $n_c = \{16, \, 64, \, 256\}$ are depicted in Figure~\ref{fig:cells2d}. One ghost layer is padded around the plate to ensure completeness near the plate's boundaries. We consider linear and quadratic local polynomials $r_p \in \{1, \, 2\}$ mollified with a $C^2$-smooth spline mollifier obtained from the tensor product of a one-dimensional spline curve
\begin{equation} \label{eq:quartSpline}
	m(x) =  
	\begin{cases} 
		\frac{35}{16h_m} \left( 1 - 12\left( \frac{x}{h_m} \right)^2 + 48 \left( \frac{x}{h_m} \right)^4 - 64 \left( \frac{x}{h_m} \right)^6 \right) & \text{  if  } |x| < \frac{h_m}{2}\\ 
		0 & \text{  if  } |x| \geq \frac{h_m}{2}\end{cases}
\end{equation}
where the mollifier width $h_m$ is obtained by averaging the total area of the domain $\Omega$ by the total number of internal cells, that is, 
\begin{equation}
	h_m = 2 \left( \frac{1}{n_c} \right)^{0.5} \, .
\end{equation}
We determine the total number of collocation points $n_z$ according to the total number of basis functions $n_b$. In particular, we require $n_z \geq n_b$ regardless of the type of point distributions used. The number of basis functions is determined by
\begin{equation} \label{eq:numBasis2dPlate}
	n_b = \left( n_c + n_g \right) \, | \vec p_i |
\end{equation}
where $n_c$ is the number of cells and $n_g$ is the number of ghost cells. The notation $|\vec p_i|$ indicates the number of monomial basis functions in each cell, which depends on the order $r_p$ and dimension $d$. For the two-dimensional case, $|\vec p_i| = 3$ applies for the linear and $|\vec p_i| = 6$ applies for the quadratic polynomial order. For example, when the number of internal cells is $n_c = 16$ and the number of ghost cells is $n_g = 20$, the total number of basis functions for the quadratic case is $n_b = (16 + 20) \cdot 6 = 216$. 

\begin{figure}[]
	\centering
	\subfloat[][Uniform \label{fig:cell16Unif}]{
		\includegraphics[scale=0.45]{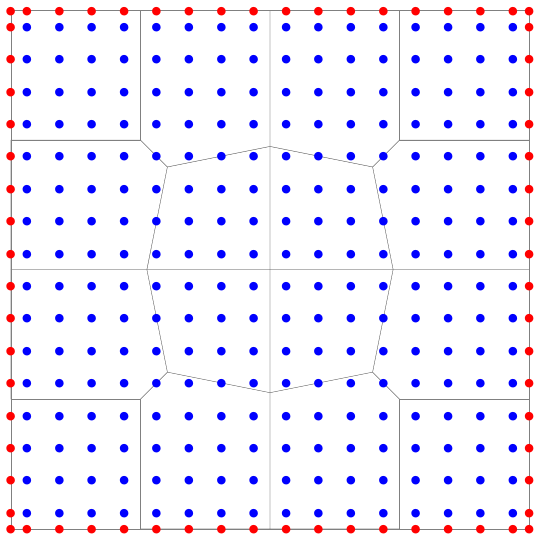}
	}
	\hspace{0.03\textwidth}
	\subfloat[][Gauss quadrature \label{fig:cell16Gauss}]{
		\includegraphics[scale=0.45]{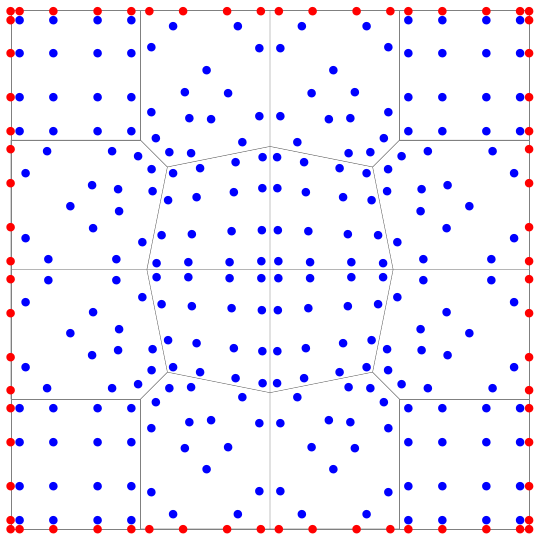}
	}
	\hspace{0.03\textwidth}
	\subfloat[][Quasi-random \label{fig:cell16Perturbed}]{
		\includegraphics[scale=0.45]{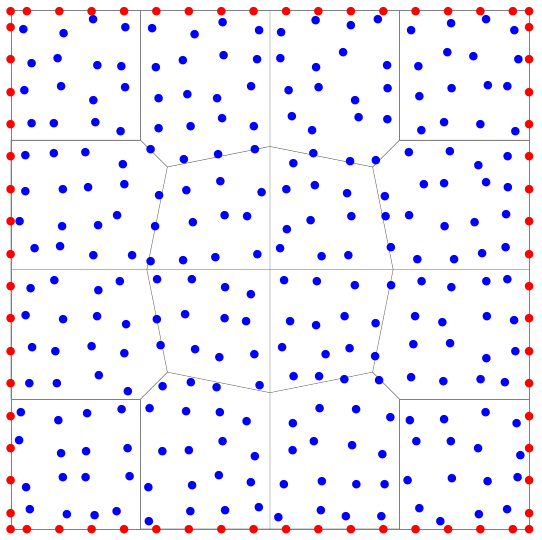}
	}
	\caption[]{Collocation points according to the uniform (a), Gauss quadrature (b), and quasi-random (c) distributions for the case with $n_c = 16$ cells.}
	\label{fig:cellGpoi}
\end{figure}

We arrange the collocation points according to the uniform, Gauss quadrature, and quasi-random distribution in a similar way as in the previous one-dimensional examples. In the uniform case, we consider a tensor product of the one-dimensional uniform point distribution over the domain $\Omega$. In the above example with $n_c = 16$ cells and the quadratic polynomial order, we select $n_z = 16^2$ collocation points as depicted in Figure~\ref{fig:cell16Unif}. For the Gauss quadrature case, we first triangulate the polytopic cells and subsequently distribute the quadrature points by mapping them from a reference triangle. In the case of quadrilateral cells, we use Gauss quadrature points mapped from a quadrilateral reference element, as illustrated in Figure~\ref{fig:cell16Gauss}. The quadrature order is uniformly chosen for all cells so that the criterion $n_z \geq n_b$ is satisfied. For the quasi-random point distribution, we first consider the uniform point arrangement and apply a small perturbation to each point. Moreover, in each coordinate axis, we consider a perturbation sampled from a uniform distribution $\epsilon \sim \set U(-\sigma, \, \sigma)$, where $\sigma$ is 15\% of the equidistant point spacing (Figure~\ref{fig:cell16Perturbed}). 

\begin{figure}
	\centering
	\includegraphics[width=6cm]{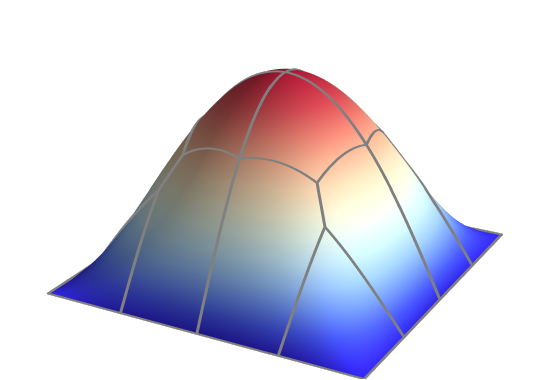}
	\caption{Solution to the elastic plate problem computed over $n_c = 16$ cell discretisation.}
	\label{fig:plateSoln}
\end{figure}

\begin{figure}[]
	\centering
	\subfloat[][$L^2$-norm error \label{fig:2d-el-L2}]{
		\includegraphics[scale=0.5]{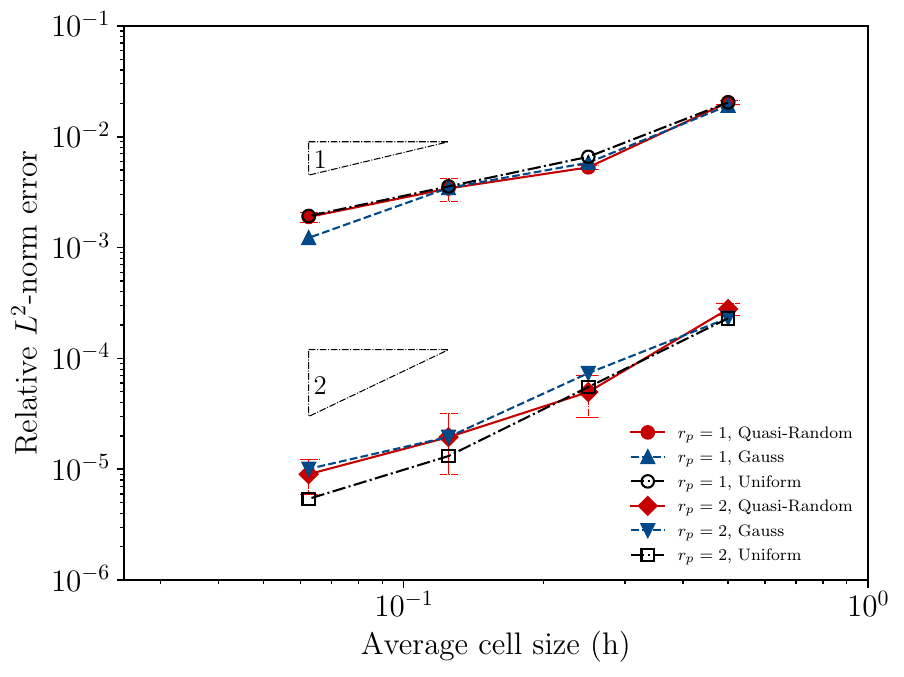}
	}
	\hspace{0.05\textwidth}
	\subfloat[][$H^1$-seminorm error \label{fig:2d-el-H1}]{
		\includegraphics[scale=0.5]{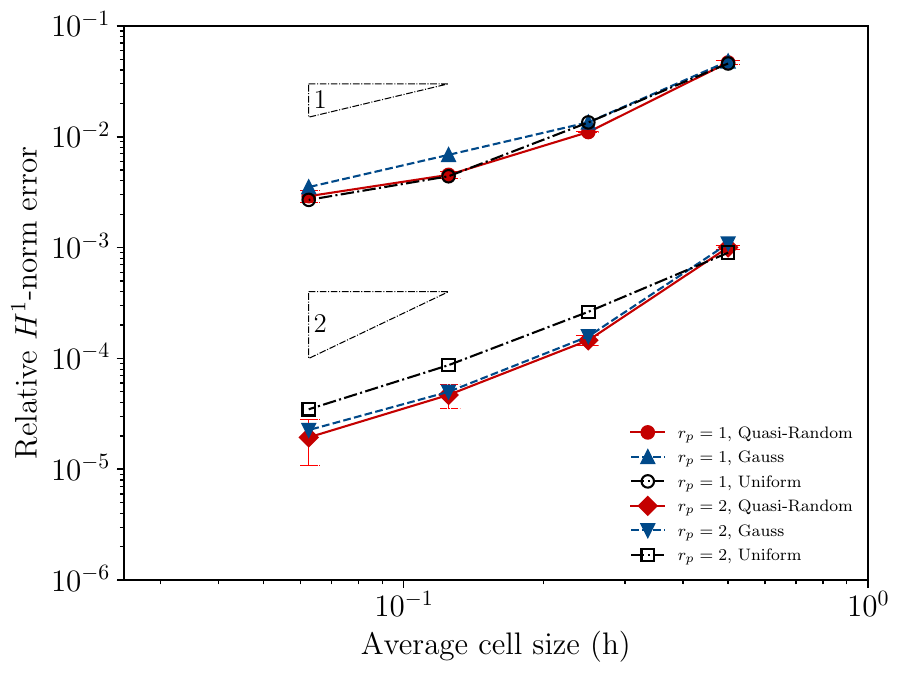}
	}
	\caption[]{Two-dimensional elastic plate problem. Convergence with a $C^2$ hexic  spline mollifier with piecewise polynomial of order $r_p = \{1, 2\}$. }
	\label{fig:plate-result}
\end{figure}

Figure~\ref{fig:plateSoln} shows the solution contour to the elastic plate problem computed over $n_c = 16$ polytopic cells. Figure~\ref{fig:plate-result} illustrates the convergence of the solution approximation in the $L^2$-norm and $H^1$-seminorm. The convergence rates of $r_p$ are achieved in both norms for both the linear $r_p = 1$ and quadratic $r_p = 2$. For the linear case $r_p = 1$, the error constants of the three point distributions are comparable, with the mean error of the quasi-random distribution having a slightly lower convergence rate.  Furthermore, for the quadratic order $r_p = 2$, the quasi-random distribution has a slightly higher mean error in the $L^2$-norm and a slightly lower convergence rate compared to the other two distributions. Unlike in the one-dimensional example, we cannot deduce a trend in the standard deviation of errors with quasi-random collocation points, which might be indirectly influenced by the non-nested refinement. 

%--------------------------------------------------------------------------------
\subsubsection{Two-dimensional plate bending} 
\label{sec:plate-bending}
%--------------------------------------------------------------------------------
%
In this section, we consider the two-dimensional bending problem on a two-dimensional square plate defined by the domain $\Omega = (0,1) \times (0,1)$ with boundary $\Gamma$. The governing biharmonic equation and the boundary conditions read:
\begin{subequations}
	\begin{align}
		D \, \Delta^2 u &= q  \qquad \text{  in  } \Omega \, , \\
		u &= 0 \qquad \text{  on  } \Gamma \, ,\\
		\nabla u \cdot  \vec n &= 0  \qquad \text{  on  } \Gamma \, .
	\end{align}
\end{subequations}
Here, $D$ is a constant associated with the material properties and the thickness of the plate. In our computations, we assume $D = 1$ and the right-hand term $q(\vec x) $ is chosen according to
\begin{equation}
	q(\vec x) = -16 \pi^4 \left( \cos \left(2 \pi x^{(1)} \right) - 4 \cos \left(2 \pi x^{(1)} \right) \cos \left(2 \pi x^{(2)} \right) + \cos \left(2 \pi x^{(2)} \right)  \right)
\end{equation}
which leads to the solution of 
\begin{equation}
	u(\vec x) = \left(1 - \cos \left( 2 \pi x^{(1)} \right) \right) \left( 1 - \cos \left(2 \pi x^{(2)} \right) \right) \, .
\end{equation}

The domain $\Omega$ is partitioned into $n_c$ non-uniform polytopic cells as used in the previous two-dimensional elastic plate cells mentioned in Section~\ref{sec:plate} and Figure~\ref{fig:cells2d}. One ghost layer is padded around the plate to ensure completeness near the plate's boundaries. We consider quartic and quintic local polynomials $r_p \in \{4, \, 5\}$ mollified with a $C^4$-smooth spline mollifier obtained from the tensor product of a one-dimensional spline curve
\begin{equation} \label{eq:decicSpline}
	m(x) =  
	\begin{cases} 
		\frac{2772}{1024 h_m} \left( 1 - 20\left( \frac{x}{h_m} \right)^2 + 160 \left( \frac{x}{h_m} \right)^4 - 640 \left( \frac{x}{h_m} \right)^6  + 1280 \left( \frac{x}{h_m} \right)^8 - 640 \left( \frac{x}{h_m} \right)^{10} \right) & \text{  if  } |x| < \frac{h_m}{2}\\ 
		0 & \text{  if  } |x| \geq \frac{h_m}{2}\end{cases} \, .
\end{equation}

The mollifier width $h_m$ is obtained by averaging the total area of the domain $\Omega$ by the total number of internal cells, like in our plate case (Section~\ref{sec:plate}). In this example, we consider a Gauss quadrature of order $\gamma$ as collocation points mapped from the simplices division of each cells. The Gauss quadrature order is set as $\gamma = 7$, which satisfies our requirement that $n_z \leq n_b$.  According to~\eqref{eq:numBasis2dPlate}, for the polynomial order $r_p = 4$, we consider the total number of monomial basis functions in one cell to be $|\vec p_i| = 15$. Likewise, for $r_p = 5$, there are $|\vec p_i| = 21$ monomials in one Voronoi cell. 

\begin{figure}
	\centering
	\includegraphics[width=6cm]{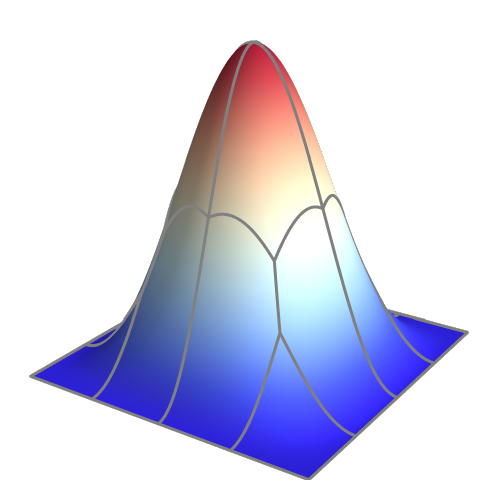}
	\caption{Solution to the two-dimensional biharmonic problem computed over $n_c = 16$ cell discretisation.}
	\label{fig:biharmonicSoln}
\end{figure}

\begin{figure}[]
	\centering
	\subfloat[][$L^2$-norm error \label{fig:2d-bi-L2}]{
		\includegraphics[scale=0.5]{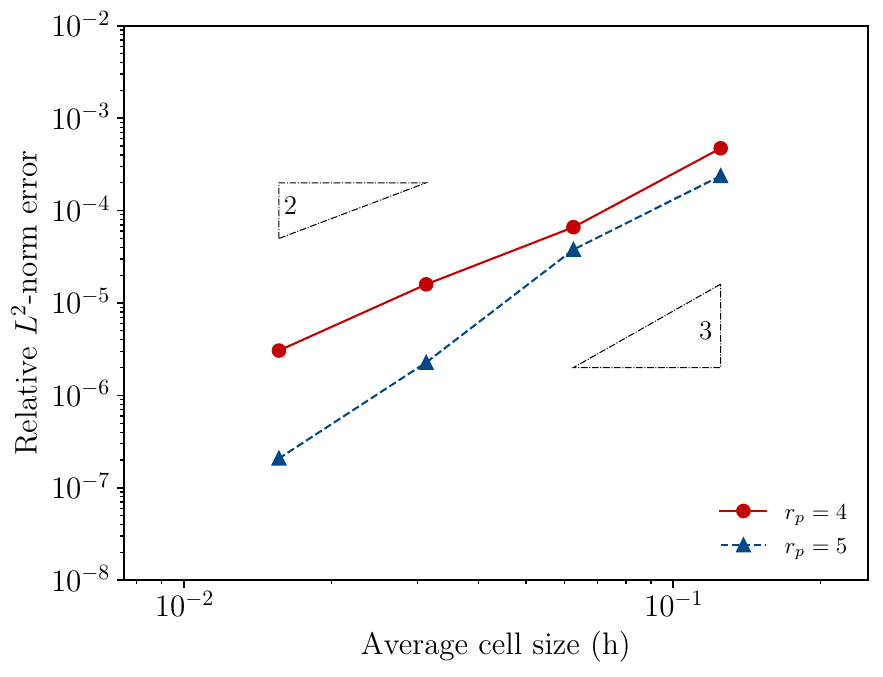}
	}
	\hspace{0.05\textwidth}
	\subfloat[][$H^1$-seminorm error \label{fig:2d-bi-H1}]{
		\includegraphics[scale=0.5]{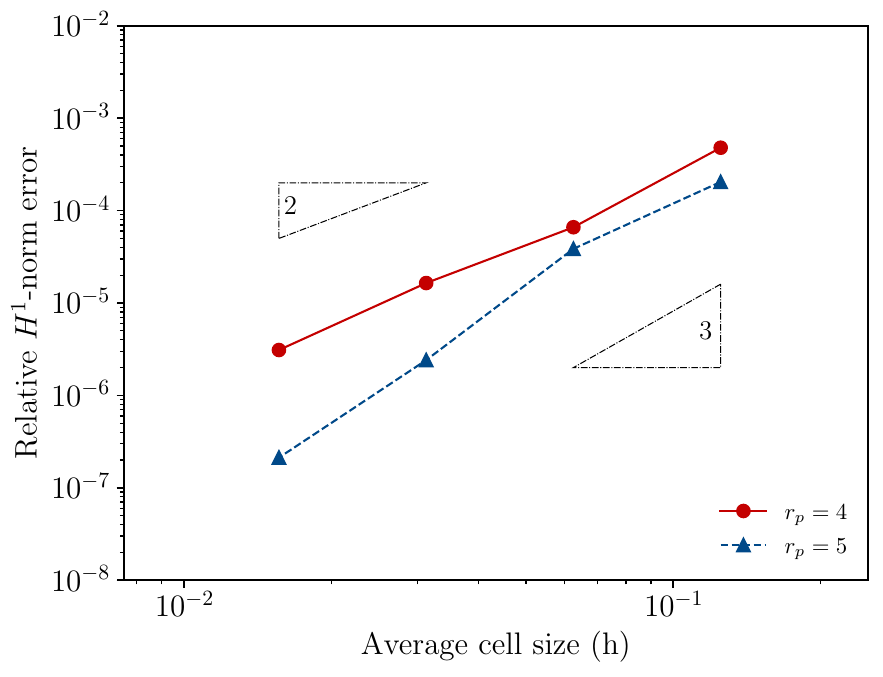}
	}
	\caption[]{Two-dimensional biharmonic plate problem. Convergence with $C^4$-smooth spline mollifier with a piecewise polynomial of order $r_p = \{4, 5\}$. }
	\label{fig:2dbiharmonic-result}
\end{figure}

Figure~\ref{fig:biharmonicSoln} shows the solution contour to the plate bending problem computed over $n_c = 16$ polytopic cells. The convergence of relative errors in the $L^2$-norm and $H^1$-seminorm is illustrated in Figure~\ref{fig:2dbiharmonic-result}. In both the $L^2$-norm and $H^1$-seminorm, the error converges approximately with a rate of $r_p - 2$ for both $r_p = 4$ and $r_p = 5$. In addition, the scaling of the monomial basis introduced in~\eqref{eq:scaledMono} is crucial for the stability of the linear system because of the high-order polynomial involved. Moreover, because of the higher order derivatives involved in this example, we use another scaling factor of $(h_m)^d$, where $d$ is the derivative order, to better condition the linear system, as explained in Section~\ref{sec:discretisation}.

%-------------------------------------------------------------------------------- 
\subsubsection{Two-dimensional infinite plate with a hole}
\label{sec:plate-hole}
\begin{figure}
	\centering
	\includegraphics[width=6cm]{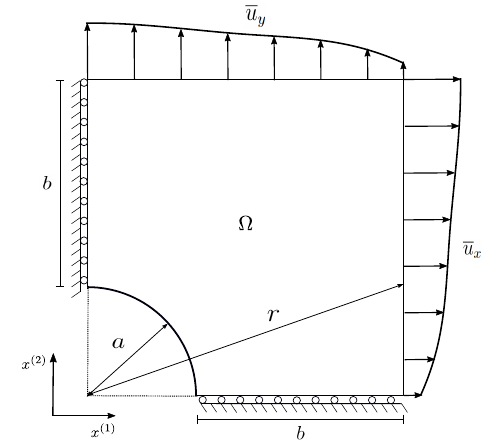}
	\caption{Schematic of the elastic plate with a hole problem.}
	\label{fig:plate-with-hole}
\end{figure}

%-------------------------------------------------------------------------------- 
%
In this section, we consider an infinite plate with a circular hole subjected to uniaxial tension. The tension $\sigma_\infty = 10^6$ is applied in the $x$-direction (Figure~\ref{fig:plate-with-hole}). Due to symmetry, we consider only a quarter of the plate with a unit length and the a hole radius of $a = 0.25$. The material has Young’s modulus of $E = 70 \times 10^6$ and its Poisson’s ratio is $\nu = 0.3$. The infinite plate with a hole problem has a closed-form analytic solution \cite{timoshenko1970} as follows
\begin{align} 
	\label{eq:plate-hole-anlyt}
	u_x &= \frac{\sigma_\infty \, a}{8 \mu} \bigg[ \frac{r}{a} (\kappa + 1) \cos (\theta) + \frac{2a}{r} \big( (1 + \kappa) \cos (\theta) + \cos (3\theta) \big) - \frac{2a^3}{r^3} \cos(3\theta) \bigg] \\
	u_y &= \frac{\sigma_\infty \, a}{8 \mu} \bigg[ \frac{r}{a} (\kappa - 3) \sin (\theta) + \frac{2a}{r} \big( (1 - \kappa) \sin (\theta) + \sin (3\theta) \big) - \frac{2a^3}{r^3} \sin(3\theta) \bigg] \, ,
\end{align}
where $r$ is the distance from the centre of the hole and $\kappa$ is the Kolosov constant
\begin{equation} 
	\label{kolosov}
	\kappa = \frac{3-\nu}{1+\nu} 
\end{equation} 
for plane stress. Dirichlet boundary conditions are imposed over the entire boundary of the plate. 

\begin{figure}[]
	\centering
	\subfloat[][$n_c = 16$ \label{fig:cell16pwh}]{
		\includegraphics[scale=0.39]{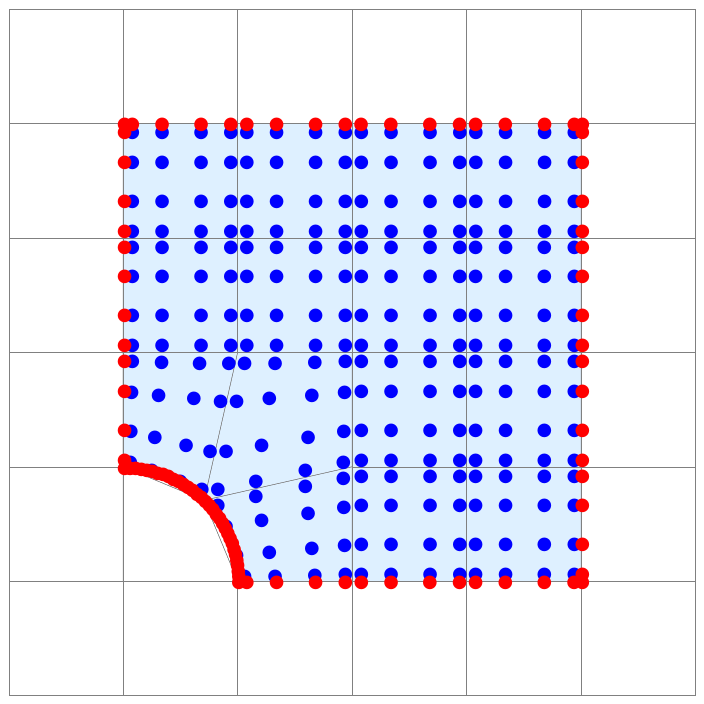}
	}
	\hspace{0.02\textwidth}
	\subfloat[][$n_c = 64$ \label{fig:cell64pwh}]{
		\includegraphics[scale=0.39]{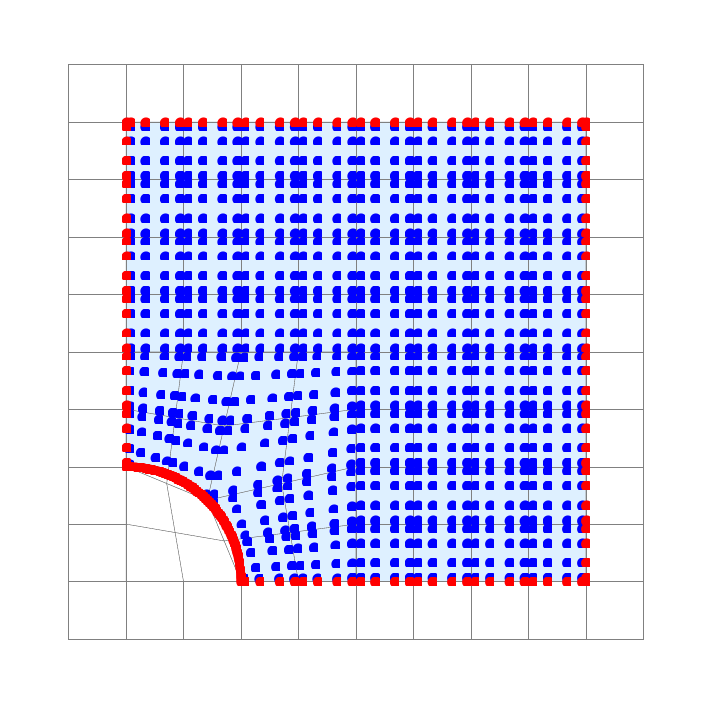}
	}
	%\hspace{0.0\textwidth}
	\subfloat[][$n_c = 256$ \label{fig:cell256pwh}]{
		\includegraphics[scale=0.39]{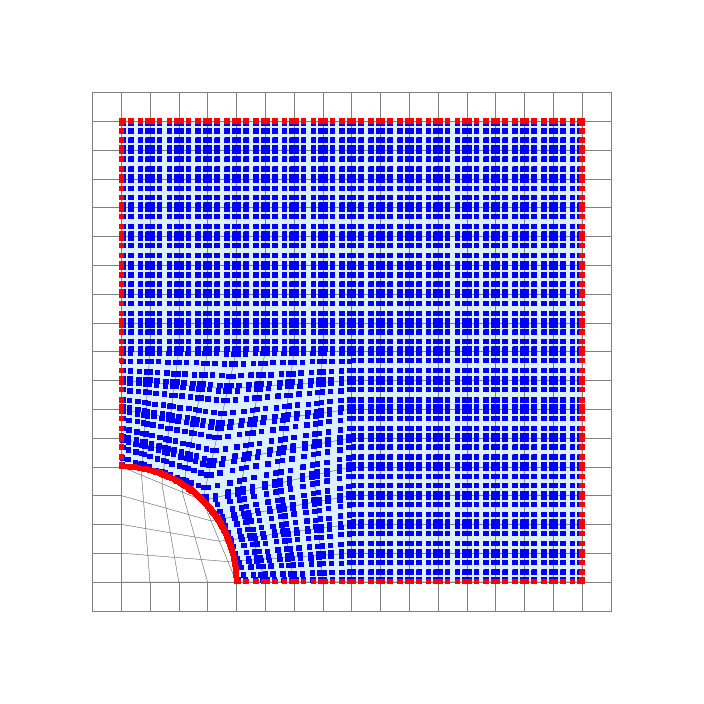}
	}
	\caption[]{Distribution of collocation points inside the plate-with-hole domain and on its boundary according to the Gauss quadrature arrangement.}
	\label{fig:voronoi-plate-hole}
\end{figure}

In this example, we consider an initial mesh of $n_c = 16$ quadrilateral cells as shown in Figure~\ref{fig:cell16pwh}. The refined meshes are obtained by introducing new vertices in the middle of each edge and subsequently subdividing the cells into four, see Figure~\ref{fig:cell64pwh} and~\ref{fig:cell256pwh} for $n_c = 64$ and $n_c = 256$, respectively. As in previous examples, ghost cells are considered to ensure completeness at the boundary and are explicitly shown in~\ref{fig:voronoi-plate-hole}. Here, we consider Gauss quadrature of order $\gamma$ as collocation points mapped from the reference quadrilateral to each cell, see Figure~\ref{fig:voronoi-plate-hole}. The Gauss quadrature order is $\gamma = 3$ and $\gamma = 4$ for linear $r_p = 1$ and quadratic $r_p = 2$ polynomial orders, respectively. Because the cell edges do not align with the domain boundary around the hole, some collocation points may lie outside of the domain $\Omega$. Such collocation points are removed from computation through an auxiliary detection method using the signed distance function $\phi(\vec x)$ of the domain $\Omega$. By definition, the signed distance function~$\phi(\vec x)$  is positive inside the domain, negative outside the domain, and the zeroth isosurface~$\phi^{-1}(0)$ corresponds to the boundary~$\Gamma$. Therefore, we consider only quadrature points that satisfy $\phi(\vec z_k) > 10^{-5}$ as internal collocation points. Furthermore, to enhance accuracy around the hole, boundary collocation points are mapped from a reference one-dimensional line according to the eight-th order Gauss quadrature $\gamma = 8$.

\begin{figure}
	\centering
	\includegraphics[scale=0.5]{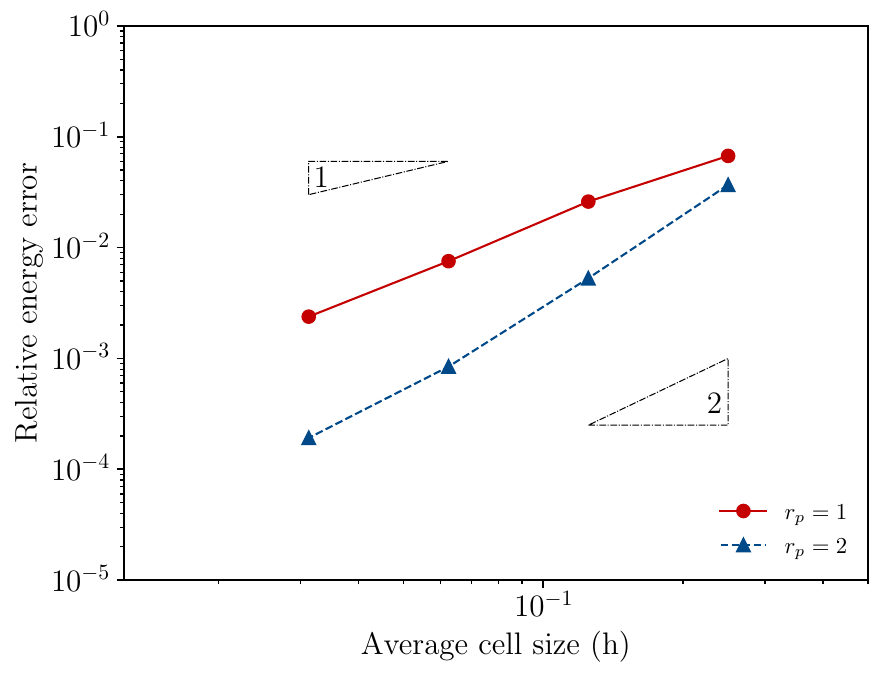}
	\caption{Elastic plate with a hole. Convergence of the relative energy norm error with a $C^2$-smooth spline mollifier and a local polynomial basis of degree $r_p = \{1, 2\}$. }
	\label{fig:energy-plate-hole}
\end{figure}

To analyse this problem, the $C^2$-smooth spline mollifier is used as described in~\eqref{eq:quartSpline}, yielding  $C^3$-smooth basis functions.  Figure~\ref{fig:energy-plate-hole} shows the convergence of the solution error in the energy norm. It is evident that convergence rate of $r_p$ is achieved for both the linear $r_p = 1$ and quadratic $r_p = 2$ polynomial orders. Furthermore, we would like to compare the basis evaluation between the proposed collocation approach and the finite element (FE) implementation~\cite{Febrianto2021}. As described in Section~\ref{sec:basis-eval}, both approaches require evaluating basis functions at Gauss quadrature points. It is important to note that the adequate quadrature order needed for the FE is determined by the polynomial order of the integrands, whereas in the collocation approach, the appropriate quadrature order is determined through a loose criterion ($n_z \geq n_b$) to avoid underdetermined matrix system. For instance, to integrate the FE stiffness term in the quadratic case $r_p = 2$ with the $C^2$ polynomial mollifier, the integrand has a maximum polynomial order of $16$. Therefore, the minimum Gauss quadrature order required is $\gamma = 9$. Assuming quadrilateral cells, leads to $81$ points per cell. By contrast, the proposed collocation approach requires $\gamma = 4$ to ensure the system is overdetermined. Moreover, although auxiliary techniques such as variationally consistent integration~\cite{Chen2013, hillman2015} can be used to aid in the FE integration, they still require an adequate base quadrature order. 

%--------------------------------------------------------------------------------  
\subsection{Three-dimensional heat transfer on a solid body}
\label{sec:solid3d}
%--------------------------------------------------------------------------------
%
%
\begin{figure}[]
	\centering
	\subfloat[][Domain definition \label{fig:boxMeshCoord}]{
		\includegraphics[scale=0.065]{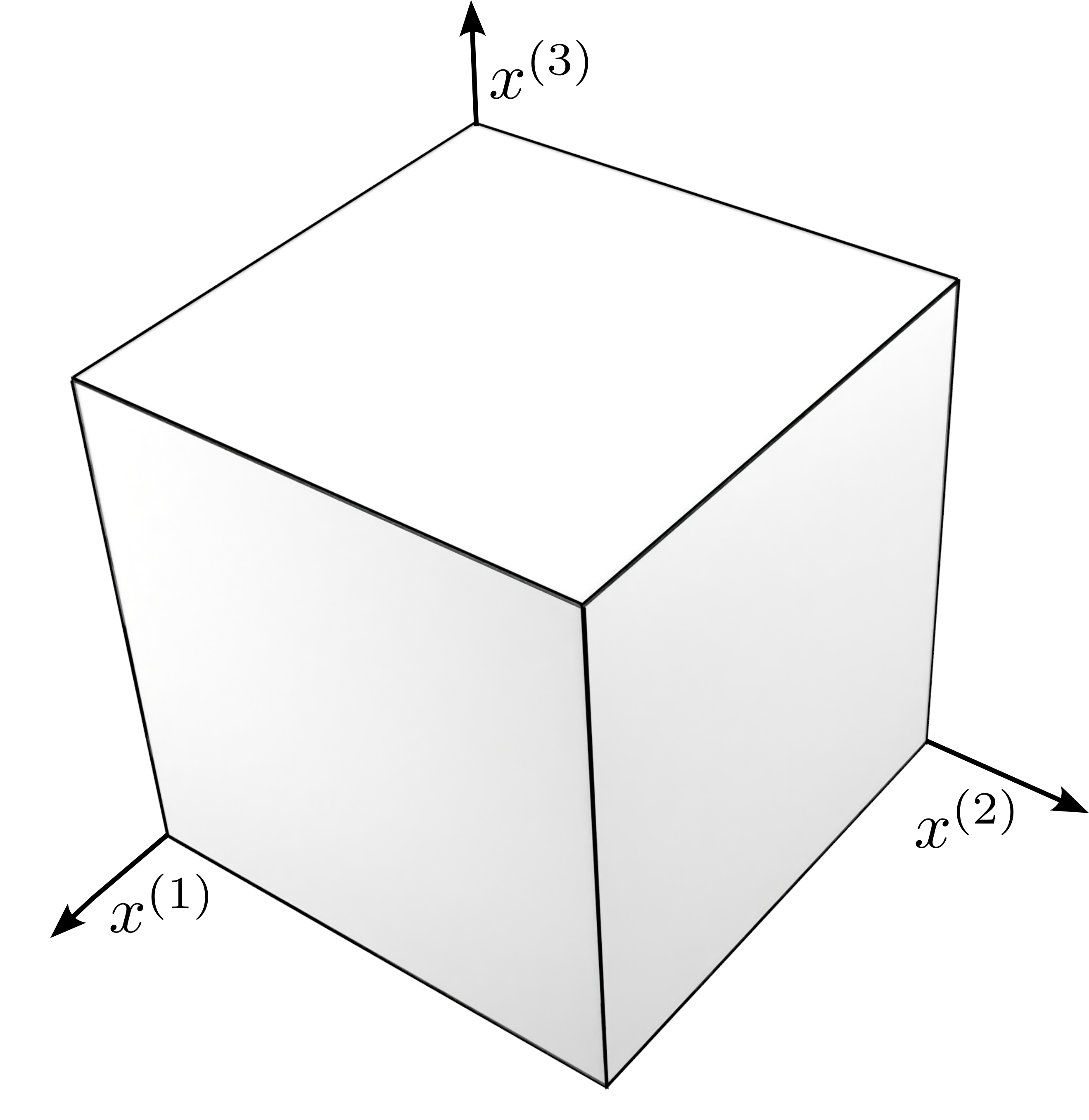}
	}
	\hspace{0.1\textwidth}
	\subfloat[][Voronoi mesh \label{fig:boxMesh}]{
		\includegraphics[scale=0.1]{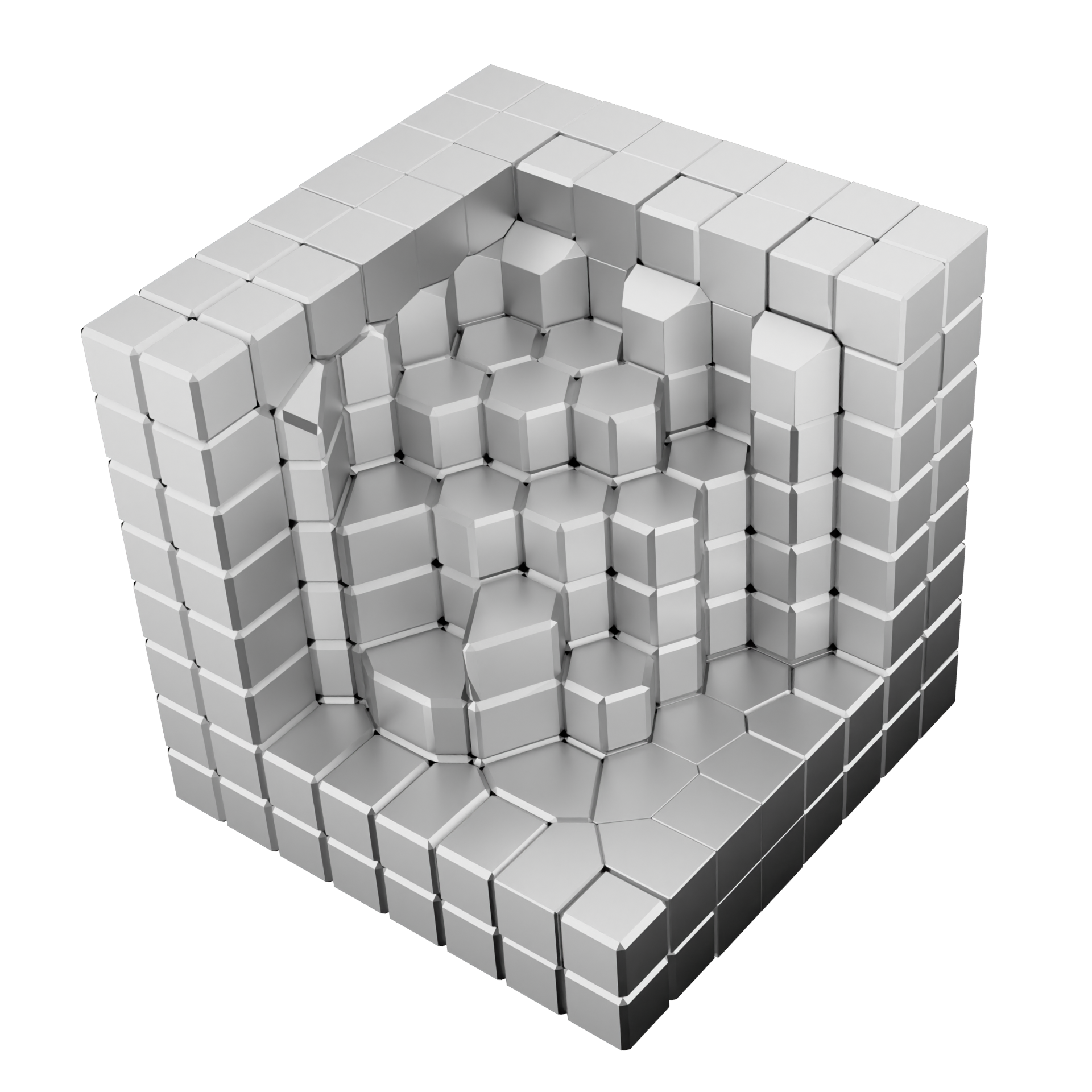}
	}
	\\
	\subfloat[][Temperature distribution over the collocation points \label{fig:boxPts}]{
		\includegraphics[scale=0.5]{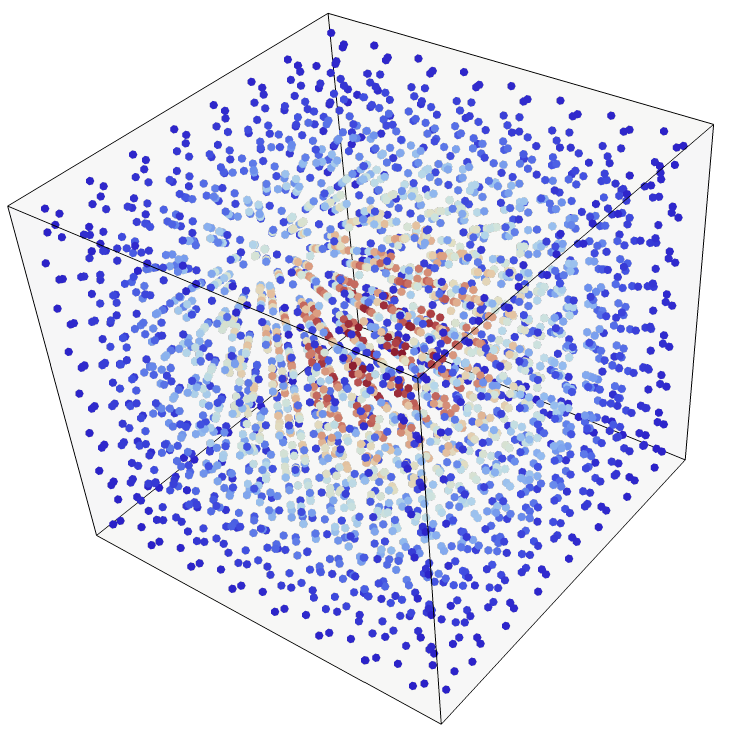}
	}
	\hspace{0.1\textwidth}
	\subfloat[][Temperature contour \label{fig:boxResult}]{
		\includegraphics[scale=0.1]{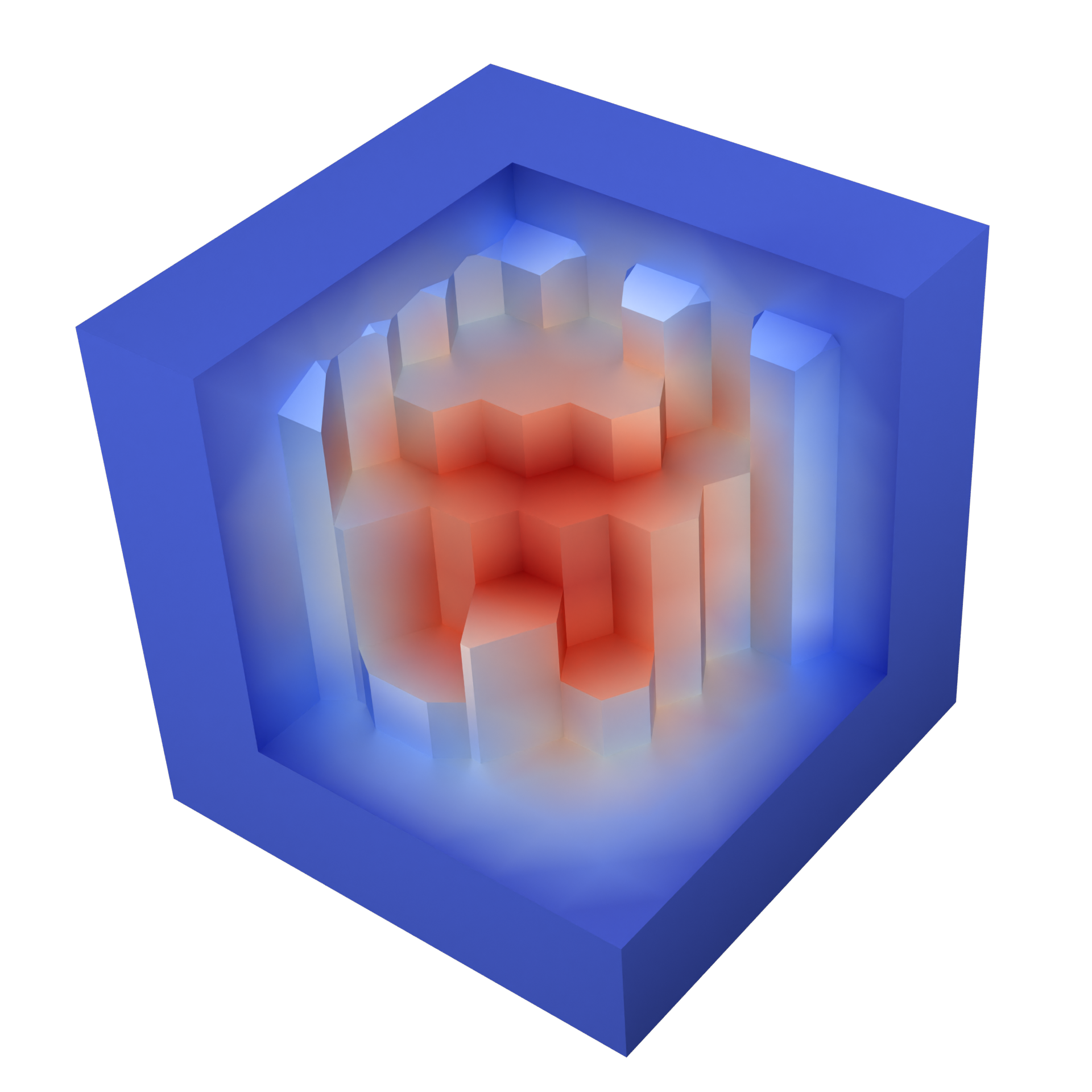}
	}
	\caption{Three-dimensional box case result. The top left image depicts the cell coordinate system, the top right image depicts the exploded Voronoi mesh that is used in this case, the bottom left image depicts the points distribution, and the bottom right image depicts the temperature contour of the result.}
	\label{fig:box}
\end{figure}

In this example, we consider a steady heat transfer with a constant coefficient of 1, that is, $\nabla^2 T (\vec x) = \vec f (\vec x)$. The source term is applied so that the temperature solution is $T(\vec x) = \sin \left( \pi x^{(1)}\right) \sin \left( \pi x^{(2)}\right) \sin \left( \pi x^{(3)}\right)$. A linear local polynomial basis function is then chosen in each cell and a tensor product of the $C^1$-smooth mollifier is used, as described in~\cite{Febrianto2021}. We use the proposed mollified-collocation method to solve the heat transfer problem over solid bodies discretised using polytopic meshes. Here we consider two objects as our domain: a unit cube and a dodecahedron with a uniform edge length of 0.3, as shown in Figure~\ref{fig:box} and Figure~\ref{fig:dod}, respectively. The objects are discretised using the  \textit{VoronoiMesh} function in \textit{Mathematica} by distributing quasi-uniform Voronoi seeds over a bounding box $\Omega_\Box$ larger than the domain $\Omega$. Subsequently, we exclude the cells corresponding to zero-valued basis in $\Omega$. The resulting active Voronoi meshes consist of 1000 cells for the box example and 622 cells for the dodecahedron. The discretisation has a total number of basis functions $n_b = 4000$ for the cube and $n_b = 2488$ for the dodecahedron. Furthermore, an auxiliary intersection algorithm is required to determine part of the Voronoi cells inside the domain $\Omega$ for distributing the collocation points. After obtaining part of the cell lying inside the domain, we first tessellate this part into tetrahedra and map Gauss points onto each tetrahedron. The boundary collocation points are obtained by mapping the Gauss quadratures onto the surface triangles. This approach results in $n_z = 8373$ collocation points for the cube example, and $n_z = 6522$ points for the dodecahedron. The isocontours of the computed temperature are shown in Figure~\ref{fig:box} and Figure~\ref{fig:dod} for the cube and dodecahedron, respectively. It is worth emphasising that our method allows for non-boundary fitting discretisation, which simplifies the domain discretisation into cells.

%% file: conclusions.tex
%--------------------------------------------------------------------------------
\section{Conclusions}
\label{sec:conclusion}
%--------------------------------------------------------------------------------
%

We presented a point collocation method that uses the smooth mollified basis functions to approximate the solutions of Poisson, linear elasticity, and biharmonic equations. The method attained high-order numerical convergence. The approximation properties of the mollified basis were characterised by the order of local polynomial approximants and the smoothness of the mollifier. The smoothness of the approximation using mollified basis functions remained intact even across meshes with irregular polytopic shapes. Here, we considered polynomial mollifiers with compact support and unit volume. To evaluate the basis functions at a point, a convolution integral is solved by first obtaining a compact integration domain, which is the intersection between the polytope and a box. They represent the support of the piecewise polynomial and the mollifier, respectively. Such a geometric intersection can be robustly computed using polytope clipping and convex hull algorithms implemented in Mathematica, Python, and similar geometry processing libraries~\cite{Barber1999, Rycroft2009, Ray2018}. In this work, we constructed an overdetermined linear system by choosing the number of collocation points to exceed the number of basis functions. Furthermore, to improve the conditioning of the system matrix, we scaled the basis functions and their derivatives. These treatments yielded good convergence properties irrespective of the mollifier type, support size, or the spatial distribution of the collocation points.  Finally, the proposed mollified collocation approach dispensed with the need for integrating the domain and boundary integrals in contrast to the mollified Galerkin approach~\cite{Febrianto2021}. 

There are several promising future extensions of the proposed mollified collocation approach. The first proposition concerns the requirement for the padded ghost cells to ensure the polynomial reproduction property of the mollified basis functions throughout the domain. Consequently, the basis functions associated with ghost cells have to be considered in the computation, which ultimately leads to an increased number of bases. A systematic approach to guarantee polynomial reproduction without needing ghost cells will improve the efficiency of the mollified-collocation method. One promising avenue includes morphing the kernel when approaching the boundary. In addition, the recent development of boundary-fitted Voronoi tessellations can be incorporated into the mollified-collocation framework. Efficient implementations of such algorithms have been reported, for example~\cite{Ray2018, Abdelkader2020}. Furthermore, an obvious extension of the proposed method is to consider local $p-$ and $h-$ refinement. Local $p-$ refinement is possible because of the  individual prescription of the local polynomial order for each cell. When Voronoi tessellation is used for domain discretisation, the $h-$ refinement becomes less obvious. One possible approach involves adding more Voronoi seeds to the area of interest and regenerating the Voronoi tessellation. The regularity of cells can then be improved using the standard Lloyd's iteration~\cite{Lloyd1982, Du1999}. Another promising future work involves exploring and establishing the connection between mollification and convolutional neural networks~\cite{LeCun1998}. This will allow for efficient uni- and multivariate basis evaluations, using open-source machine learning tools such as PyTorch and TensorFlow. Finally, the choice of local monomial basis requires appropriate scaling factors to better condition the system matrix. Some studies are available on the alternative basis functions for polytopic elements, for example~\cite{Sukumar2006, Schneider2019, Bunge2022}, which could be adapted for application in the mollified context. 